\documentclass[10pt]{elsarticle}

\nonstopmode

\usepackage{subcaption}

\usepackage{microtype}
\usepackage{booktabs}

\usepackage[utf8]{inputenc}
\usepackage{graphicx}
\graphicspath{{paper_im/}}
\usepackage{color}
\usepackage{verbatim}
\usepackage{chngpage}
\usepackage{hyperref}
\usepackage{url}

\usepackage{hypernat}
\usepackage{mathpazo}
\usepackage[table]{xcolor}
\usepackage{pgfplots}
\usepackage{amsmath}
\usepackage{cleveref}
\usepackage{amsfonts}
\usepackage{amssymb}

\numberwithin{equation}{section}

\usepackage{setspace}

\usepackage{stmaryrd}

\usepackage{tikz}
\usetikzlibrary{positioning}
\usetikzlibrary{calc}
\usetikzlibrary{chains}
\usetikzlibrary{shapes}
\usetikzlibrary{arrows}

\pgfdeclarelayer{background}
\pgfdeclarelayer{foreground}
\pgfsetlayers{background,main,foreground}

\usepackage{floatrow}
\newfloatcommand{capbtabbox}{table}[][\FBwidth]

\usepackage[noend]{algorithmic}
\usepackage{algorithm}
\algsetup{indent=2em}

\usepackage{enumitem}

\usepackage{multirow,hhline}

\usepackage[margin=10pt,font=small,labelfont=bf,labelsep=period]{caption}

\def\revanswera#1{#1}
\def\revanswerb#1{#1}
\def\revanswerbb#1{#1}

\begin{document}

\title{  Multi-Rate Time Integration on Overset Meshes
}

\author[uiuc-aero]{Cory Mikida}\ead{cmikida2@illinois.edu}
\author[uiuc-cs]{Andreas Kl\"ockner}\ead{andreask@illinois.edu}
\author[uiuc-aero]{Daniel Bodony}\ead{bodony@illinois.edu}
\address[uiuc-aero]{Department of Aerospace Engineering,
  University of Illinois at Urbana-Champaign, 104 South Wright St,
  Urbana, IL 61801}
\address[uiuc-cs]{Department of Computer Science,
  University of Illinois at Urbana-Champaign, 201 North Goodwin Ave,
  Urbana, IL 61801}

\begin{abstract}
Overset meshes are an effective tool for the computational fluid dynamic
simulation of problems with complex geometries or multiscale spatio-temporal
features.  When the maximum allowable timestep on one or more meshes is
significantly smaller than on the remaining meshes, standard explicit time
integrators impose inefficiencies for time-accurate calculations by requiring
that all meshes advance with the smallest timestep.  With the targeted use of
multi-rate time integrators, separate meshes can be time-marched at independent
rates to avoid wasteful computation while maintaining accuracy and stability.
This work applies time-explicit multi-rate integrators to the simulation of the
compressible Navier-Stokes equations discretized on overset meshes using
summation-by-parts (SBP) operators and simultaneous approximation term (SAT)
boundary conditions.  We introduce a class of multi-rate Adams-Bashforth
(MRAB) schemes that offer significant stability improvements and computational
efficiencies for SBP-SAT methods.  We present numerical results that confirm
the efficacy of MRAB integrators, outline a number of
implementation challenges, and demonstrate a reduction in computational cost
enabled by MRAB.  We also investigate the use of our method in the setting of a
large-scale distributed-memory parallel implementation where we discuss
concerns involving load balancing and communication efficiency.
\end{abstract}


\maketitle

\onehalfspacing

\section{Introduction} \label{intro}

In the time-explicit direct numerical simulation (DNS) of computational fluid
dynamic problems, the maximum timestep allowable for stable integration of the
governing equations is often limited by the well-known Courant-Friedrichs-Lewy
number (CFL),
\begin{align}
\text{CFL} = \frac{c\Delta t}{\Delta x} \longrightarrow \Delta t \leq \frac{\text{CFL} \cdot \Delta x}{c} \notag
\end{align}
where $\Delta x$ is the minimum grid spacing, and $c$ is the maximum characteristic speed of the physical phenomena being simulated, or its viscous analog, 
\begin{align}
\text{CFL}_{\text{visc}} = \frac{\nu \Delta t}{\Delta x^{2}} \longrightarrow \Delta t \leq \frac{\text{CFL}_{\text{visc}} \cdot \Delta x^{2}}{\nu} \notag
\end{align}
where $\nu$ is the characteristic diffusivity.  The timestep taken when
integrating over the entire computational domain using a standard explicit
single-rate integrator can be limited by dynamics that occur on a
small portion of the domain, be it due to shorter time scales of certain
physical phenomena, such as high-speed flow, or locally-high grid resolution.

With multi-rate integration, groups of degrees of freedom, or even individual
right-hand side contributions corresponding to evolution on small timescales,
can be integrated with independent time steps, allowing computational work to
be avoided on the ``slow'' components while the ``fast'' components remain
stable and well-resolved in time.  This has the potential to improve the
efficiency of the application.  The present study focuses on an application of
multi-rate integrators to problems with overset meshes.

The method of overset meshes (also known as the Chimera
method~\citep{steger1991chimera}) is an approach that discretizes a domain
using multiple independent, overlapping structured meshes, each of which is
associated with a smooth but potentially geometrically complex mapping.  The
earliest known appearance of a composite mesh method was to solve elliptic
partial differential equations~\citep{volkov1968method}, but similar methods
were soon introduced for inviscid transonic flow~\citep{magnus1970inviscid} and
the Euler equations of gas dynamics~\citep{benek1983flexible}.
\revanswera{Some of the earliest analysis of stability was performed by
Starius, who examined the stability of composite mesh methods for both elliptic
boundary value problems~\citep{starius1977composite} and hyperbolic systems of
equations~\citep{starius1980composite}.}   In the years since, high-order
overset-grid (HO-OG) approaches have been
developed\revanswera{~\citep{sherer2005high, chesshire1990composite,
chesshire1994scheme}} and applied to numerous
problems\revanswera{~\citep{lee2002high, sherer2004implicit,
chou2003compact, henshaw2006high, henshaw2009composite, sjogreen2007variable}.}
The individual meshes are structured, so that the entire discretization can be
considered locally structured and globally unstructured.  A concise summary of
recent developments in the usage of overset meshes to solve problems involving
compressible viscous fluids, along with a discussion of stable and accurate
interpolation, is given in~\citep{bodony2011provably}.  Additionally, a
discussion of various applications of overset meshes can be found
in~\citep{noack2005summary}.

\revanswera{Adaptive mesh refinement methods are similar to overset grid
methods in that they address the need for locally high spatial resolution while
also inherently introducing timescale disparity~\citep{berger1984adaptive}.  In
particular, Berger et.~al applied an automatic adaptive mesh refinement method
to a hydrodynamic shock problem in two dimensions~\citep{berger1989local},
along with hyperbolic conservation laws in three
dimensions~\citep{bell1994three}.}

Some of the earliest work on multi-rate multistep methods was done by
Gear~\citep{gear1974multirate}.  Gear and Wells~\citep{gear1984multirate} later
worked with these methods to automate the choice of step size and the partition
of fast and slow components.  The primary conclusion therein was that the
problem of automating these integrator characteristics dynamically was mostly a
function of software organization.  Engstler and
Lubich~\citep{engstler1997multirate} used Richardson extrapolation along with
simple Euler schemes to create a dynamically partitioned process of multi-rate
integration.

Local timestepping methods are related to multi-rate time integration.
Osher and Sanders~\citep{osher1983numerical} introduced numerical
approximations to conservation laws that changed the global CFL restrictions to
local ones.  Dawson and Kirby~\citep{dawson2001high} performed local time
stepping based on the work of Osher and Sanders, attaining only
first-order accuracy in time.  Tang and Warnecke~\citep{tang2006high} expanded
on this work to produce second-order accuracy via more refined projections of
solution increments at each local timestep.

G\"unther, Kv\ae rn\o, and Rentrop~\citep{gunther2001multirate} introduced
multi-rate partitioned Runge-Kutta (MPRK) schemes, starting from a discussion
of Rosenbrock-Wanner methods, and based on strategies introduced by G\"unther and
Rentrop~\citep{gunther1993multirate}. These methods focus on coupling the fast
and slow solution components primarily via interpolation and extrapolation of
state variables, echoing the earlier work in~\citep{gear1974multirate}
and~\citep{gear1984multirate}.

Savcenco et.~al~\citep{savcenco2007multirate} developed a self-adjusting
multi-rate timestepping strategy primarily using implicit Rosenbrock methods on
stiff ODEs.  Constantinescu and Sandu~\citep{constantinescu2007multirate}
developed multi-rate timestepping methods for hyperbolic conservation laws
based on a method of lines (MOL) approach with partitioned Runge-Kutta schemes.
More recently, Sandu and Constantinescu~\citep{sandu2009multirate} developed
explicit multi-rate Adams-Bashforth methods similar to the ones discussed in
this work, and while the resulting methods are exactly conservative, the
computational experiments shown are limited to one-dimensional hyperbolic
conservation laws, and the temporal accuracy of the methods is limited to
second order by an interface region.

Recently, Seny and Lambrechts extended the explicit multi-rate Runge-Kutta
schemes of Constantinescu and Sandu~\citep{constantinescu2007multirate, sandu2009multirate}
to apply to discontinuous Galerkin computations for large-scale geophysical flows, introducing the
method~\citep{seny2010multirate} and subsequently demonstrating its efficacy in
both serial and parallel implementations.  Their latest
work~\citep{seny2014efficient} focuses on efficient parallelization of the
method using a multi-constraint partitioning library.

While all of these works characterize the stability of the multi-rate
integrators numerically, only a few of them approach the topic of stability
theoretically.  An analytical discussion of the stability of a multi-rate
method (namely, for numerical integration of a system of first-order ODEs that
can be readily separated into subsystems) is given by Andrus
in~\citep{andrus1993stability} whose study combines
a fourth-order Runge-Kutta scheme (for the fast component) with a similar
third-order scheme (for the slow component), introduced by the same author
in~\citep{andrus1979numerical}.  Kv\ae rn\o also derived analytical expressions to
find stability regions for her multi-rate Runge-Kutta schemes
in~\citep{kvaerno2000stability}.

In the present study, we use a multi-rate Adams-Bashforth integrator with
enhanced stability in a distributed-memory parallel fluid solver, taking
advantage of the solver's overset mesh capabilities to segregate solution
components with differing timescales.  We do so with particular focus on the
resulting improvement in performance, reduction in work through fewer
right-hand-side evaluations, and accompanying changes in stability and
accuracy.  What results from this effort is a number of conclusions about the
viability and extensibility of these integrators to other problems and
applications.  \revanswera{We couple our multi-rate Adams-Bashforth integrators
	with a summation-by-parts--simultaneous-approximation-term (SBP-SAT)
	method~\citep{strand1994summation, carpenter1993time,
	mattsson2004stable, svard2007stable, svard2008stable,
bodony2010accuracy}, in particular with the newly-developed SBP-SAT-based
interpolation method~\citep{sharan2018time}, to produce an accurate, efficient,
and high-order approach to solving PDE-based problems.  Our approach
generalizes prior work on non-overset-based multi-rate
methods~\citep{dawson2001high, tang2006high, osher1983numerical,
constantinescu2007multirate, sandu2009multirate} and removes constraints on
temporal order of accuracy.}  We also introduce extended-history
Adams-Bashforth schemes as a basis for a new class of multi-rate integrators
that demonstrate improved stability compared to standard multi-rate
Adams-Bashforth (MRAB) schemes.  For further reference on some of the multi-rate
Adams-Bashforth schemes implemented here, thorough analyses (especially
empirical analyses of stability and the effect of various method design
choices) are given in Kl\"ockner~\citep{klockner2010high}.  Some related methods are also
discussed in a particle-in-cell context by Stock~\citep{stock2009development}.

\section{Background} \label{background}

\subsection{The Compressible Navier-Stokes Equations on Overset Meshes} \label{ns}

We apply multi-rate integrators to the compressible Navier-Stokes equations on
overset meshes, making use of an SBP-SAT discretization with a SAT-based method
that applies the effect of an interpolation as a penalty term in the right-hand
side.  Details on the Navier-Stokes formulation, including the governing
equations, nondimensionalization, and coordinate mappings can be found
in~\citep{sharan2016time}.  Details on SBP operators can be found
in~\citep{strand1994summation, carpenter1993time, mattsson2004stable}.  For the
results given later, we use the third-order SBP operator given
in~\citep{mattsson2004stable}.  Details on the accompanying SAT boundary
conditions can be found in~\citep{svard2007stable, svard2008stable,
bodony2010accuracy}.  \revanswerb{For single-grid realizations of the
scheme, the SBP-SAT discretizations we use have been shown~\citep{svard2007stable,
	svard2008stable} to be provably stable for the compressible
Navier-Stokes equations in three dimensions.}  \revanswerbb{For overset-grid
realizations in two or three dimensions, no energy stability proof has yet
been found.}

\subsection{SAT-Based Interpolation} \label{interp}

\revanswera{The method used for interpolation between overset meshes is
critical to the use of multi-rate integrators, as it forms the coupling between
the fast and slow components.  For consistency, it is important that communication
between grids occurs through right-hand side increments rather than by direct
interpolation and injection of state, given that the
multi-rate schemes we will develop later are based upon multi-step schemes that
depend directly upon accurate right-hand side histories and strictly ODE-based
time advancement of the state.}  

\begin{figure}[h!]
\centering
  \includegraphics[width=0.95\linewidth]{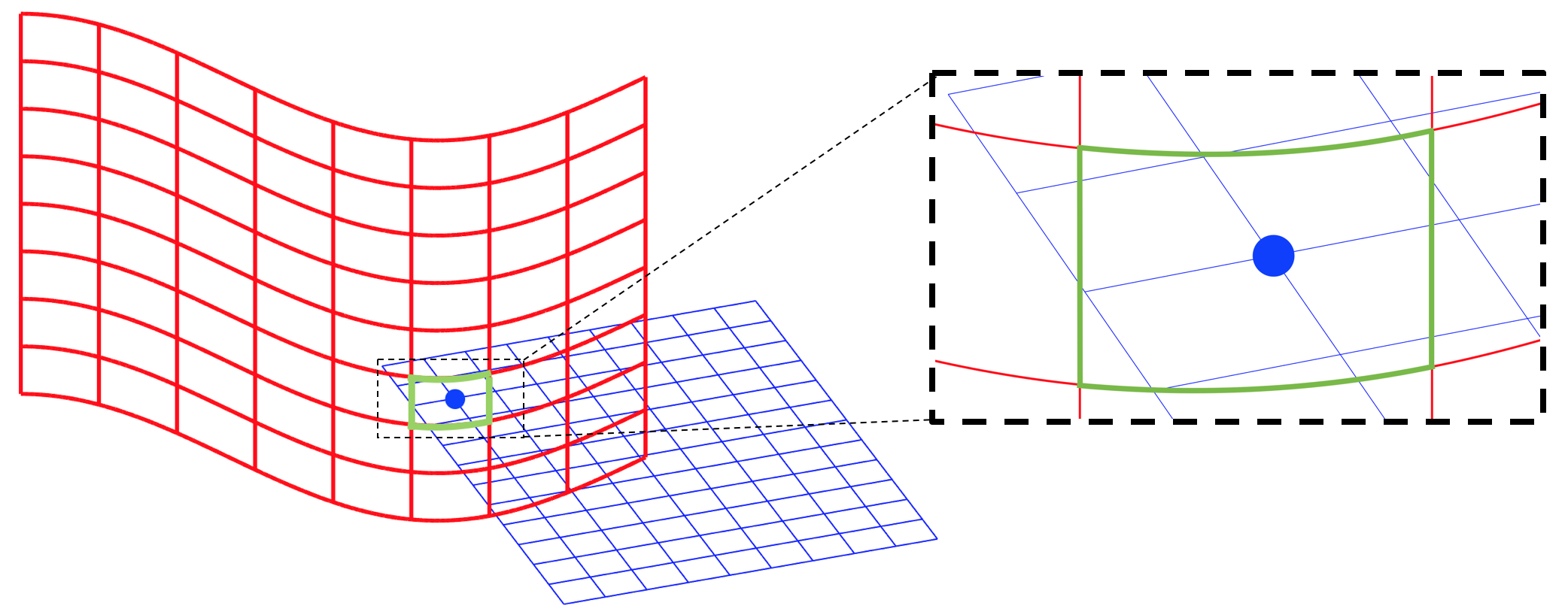}
  \caption{Sample grids for explanation of overset interpolation.  The cell outlined in green and the blue point represent an example donor cell receiver point pair.}
\label{fig:overset_new}
\end{figure}

Our interpolation relies on a Chimera framework inspired by the approach taken by the
PEGASUS~\citep{suhs2002pegasus} and BELLERO~\citep{sherer2006automated} tools.  
In general, the process of communication between grids can be broken down into
a number of phases:
\begin{itemize}
\item \textbf{Establishing communication between grids.}  Compute the bounding box for a given grid, and use it to determine collisions with any other grids.
\item \textbf{Hole cutting/fringe determination.}  Identify points on each grid
	as "fringe points" which will donate and receive data from other grids,
	and identify points on coarser background grids which are well within
	the boundaries of the finer feature grids, and can thus be deemed
	inactive.
\item \textbf{Donor-receiver pair search --- see Figure~\ref{fig:overset_new}.}  Fringe points on the receiver grid are paired with donor cells on the donor grid.
\item \textbf{Interpolation.}  State data from the points in the donor cell is
	transferred to the receiver point via Lagrangian interpolation, with
	corresponding weights determined as a function of Lagrange shape
	functions.
\end{itemize}

See~\citep{bodony2011provably} for further discussion of this implementation,
including algorithms for hole-cutting and the donor-receiver pair search.  Once
the donor data is sent (interpolated) to the receiver grids, it must be applied
to the receiver state.  The conventional method of interpolation between
overset meshes uses an ``injection'' procedure to apply the result of overset
interpolation, in which the state values on the receiver grid are directly
overwritten by the interpolated values from the donor grid.  The interpolation
itself is performed at each right-hand side evaluation.  \revanswerb{This
``injection" procedure is demonstrated in~\citep{sharan2018time} to be
unstable for the Euler equations in two dimensions.}
\revanswera{Injection-based interpolation schemes have been also used with
SBP operators by Sj\"ogreen et.~al~\citep{sjogreen2007variable} and with
adaptive mesh refinement codes~\citep{berger1984adaptive, berger1989local,
bell1994three} and have obtained conservation via interpolation of fluxes
rather than state.}

An alternative interface treatment given in~\citep{sharan2018time} follows a
methodology similar to the weak boundary treatment using
SAT~\citep{svard2007stable, svard2008stable, bodony2010accuracy} and applies
the interpolated values as a penalization term in the right-hand side via a
target vector.  Using the same notation as~\citep{sharan2018time}, we describe
the overlapping interface as a $\kappa^{\pm}$ boundary where 
$\kappa = \xi$, $\eta$, or $\zeta$.  $\kappa$ is the normal
direction to the face on which the grid point lies, and the $\pm$ superscript
indicates inflow ($+$) or outflow ($-$).  If $\mathbf{q}_{ijk}$ denotes the
solution at a grid point on the interface, with the interpolated value
from the donor grid given as $\hat{\mathbf{q}}_{ijk}$, we can express the
discretized equation at this point as
\begin{align}
\frac{d\mathbf{q}_{ijk}}{dt} = -(D_{\xi m}\mathbf{F}_{m})_{ijk} -& p_{0}^{-1}(\sigma^{I}K_{\kappa}^{\pm} + \sigma_{1}^{V}I_{5})(\mathbf{q}_{ijk} - \hat{\mathbf{q}}_{ijk}) \notag \\ 
+& \sigma_{2}^{\pm V}\left((F_{\kappa}^{V})_{ijk} - (\hat{F}_{\kappa}^{V})_{ijk}\right), \label{eq:int}
\end{align}
where $(D_{\xi_{m}} \mathbf{F}_{m})_{ijk}$ denotes the derivatives of the
fluxes, $\mathbf{F}_{m} = \mathbf{F}_{m}^{I} - \mathbf{F}_{m}^{V}$, $p_{0}$ is
the (1,1) element of the positive-definite matrix associated with the SBP
operator, $I_{5}$ is an identity matrix of size $5 \times 5$, and
$K_{\kappa}^{\pm} = T_{\kappa}\left(\frac{|\Lambda_{\kappa}| \pm
\Lambda_{\kappa}}{2}\right)T_{\kappa}^{-1}$, where $T$ and $\Lambda$ are the
transformation and diagonal matrices defined in~\citep{svard2007stable,
svard2008stable, bodony2010accuracy}.  $(F_{\kappa}^{V})_{ijk}$ denotes the
viscous flux at the interface point, and all "hatted" terms indicate
interpolated values.  Note also that if the grid point lies on an edge or a
corner in 3 dimensions, the interface terms for \emph{each} of the three directions
normal to the edge, or all directions in the case of a corner, must be added.
The penalty parameters $\sigma$ used in \eqref{eq:int} are chosen as
\begin{align}
\sigma^{I} = \frac{1}{2}, \hspace{4mm} \sigma_{1}^{V} = \frac{1}{2\text{Re}}(\kappa_{x}^{2} + \kappa_{y}^{2} + \kappa_{z}^{2}), \hspace{4mm} \sigma_{2}^{\pm V} = \pm \frac{1}{2} \notag
\end{align}
for an inflow ($+$) or outflow ($-$) interface point.  \revanswera{This
scheme is shown in \citep{sharan2018time} to be provably time-stable for
\revanswerbb{one-dimensional hyperbolic systems using} SBP-SAT discretizations
of up to global fourth-order accuracy on overlapping grids.  In addition, a
number of numerical examples, including
extension of the scheme to both the two-dimensional Euler equations and the
three-dimensional compressible Navier-Stokes equations for both Cartesian and
curvilinear grids, demonstrate both stability and superior numerical
characteristics compared to the injection method while also removing the need
for artificial damping or filtering.}  \revanswerbb{Furthermore, proper
convergence to the design order of accuracy is demonstrated with the use of
this scheme in conjunction with SBP 1-2-1, SBP 2-4-2, and SBP 3-6-3 schemes on
the Euler equations in multiple dimensions.}

\revanswera{The injection-based method described
in~\citep{bodony2011provably}} is incompatible with the strictly ODE-based time
advancement required for the schemes we will describe.  Multi-step schemes are
explicitly reliant on maintaining accurate right-hand side histories, whereas
injection methods rely on in-place modification of state variables.  The SAT
interface treatment of \citep{sharan2018time} instead applies the effects of
inter-grid communication as a right-hand side penalty, and is the method we
will use.

\revanswerbb{The expected order of convergence for overset cases using this
interpolation scheme in conjunction with a 3rd-order SBP-SAT spatial discretization is demonstrated using the method of manufactured solutions for the right-hand-side of the compressible Navier-Stokes equations (Eq.~\ref{eq:int}) as well as on analytical solutions for the inviscid (convecting vortex) and viscous (shock wave) forms of the equations.  These results are given in \ref{convergence}.}

\section{Time Integration} \label{integration}

\subsection{Adams-Bashforth Integration} \label{ab}

To fix notation for new schemes developed later, we give a brief
derivation of a standard Adams-Bashforth (AB) integrator, as described
in~\cite{bashforth1883attempt}.  We start with a model IVP given by 
\begin{align}
\frac{\text{d}y}{\text{d}t} = F(t,y), \quad y(0) = y_{0}. \notag
\end{align}
This is the form that results from a method of lines (MOL) approach to solving
time-dependent partial differential equations like those considered in this
study.  We approximate the time dependency of the right-hand side function with
a polynomial with coefficients $\boldsymbol{\alpha}$ (formed by interpolating
past values of $F(t,y)$), extrapolate with that polynomial approximation, and
integrate the extrapolant.  We use a Vandermonde matrix to construct a linear
system to obtain the coefficients $\boldsymbol{\alpha}$ to be used in
extrapolation from history values:
\begin{align}
V^{T} \cdot \boldsymbol{\alpha} = \int_0^{\Delta t} \tau^{i} d\tau = \frac{(\Delta t)^{i}}{i+1}, \quad i = 1,2...,n, \quad V = \begin{bmatrix} 
    1 & t_{1} &  \hdots   & t_{1}^{n-1}  \\
    1 & t_{2} & \hdots &  t_{2}^{n-1} \\
      \vdots  & \vdots  &  \ddots   &  \vdots   \\
       1 &   t_{n}  &  \hdots  & t_{n}^{n-1} \\   
        \end{bmatrix}, \label{eq:vandermonde}
\end{align}
where $\int_0^{\Delta t} \tau^{i} d\tau$ is a vector evaluating the integral of
the interpolation polynomial, and $V$ is the Vandermonde matrix with monomial
basis and nodes $t_{1}, t_{2}, \hdots t_{n}$, corresponding to past time
values. In \eqref{eq:vandermonde}, $n$ is equal to the order of the integrator,
and $t_{i}$ are the time history values, with $0 \leq t_{1} < t_{2} \hdots <
t_{n}$.  The coefficients $\boldsymbol{\alpha}$ are used to extrapolate to the
next state via
\begin{align}
y(t_{i+1}) = y(t_{i}) +& \alpha_{1}F(t_{i-n+1},y_{i-n+1}) + \alpha_{2}F(t_{i-n+2},y_{i-n+2}) \notag \\
+& \cdots + \alpha_{n}F(t_{i},y_{i}). \label{eq:ab_general}
\end{align}
Clearly, the length of the past history needed to calculate a step (and, thus, the memory required) influences the order of accuracy attained.

An alternative time integration method is required for the first few time steps
(the exact number of which is dependent on the number of history values needed)
in order to establish right-hand side history and "bootstrap" the method.  We
use a third-order Runge-Kutta (RK3) integrator~\cite{heun1900neue} to bootstrap
the third-order AB methods, whereas a fourth-order Runge-Kutta (RK4)
integrator~\cite{kutta1901beitrag} is used to bootstrap the fourth-order AB
methods.

\subsection{Extended-History Adams-Bashforth Schemes} \label{ab_new}
We will see in Section~\ref{sat_stab} that the ODE systems resulting from our
SBP-SAT discretization on overset meshes yield eigenvalue spectra that extend
far along the negative real axis.  To improve the stability of the AB schemes
for ODE systems containing such eigenvalues, we develop new
\emph{extended-history Adams-Bashforth schemes}.  As shown in
\eqref{eq:vandermonde}, the standard $n$-order AB scheme constructs a square
($n \times n$) Vandermonde matrix to produce a linear system that, when solved,
gives the vector of $n$ coefficients $\boldsymbol{\alpha}$ to be used in
extrapolation from history values.  With the extended-history scheme, we modify
a family of AB schemes to include $m$ history values with $m>n$, resulting in a
non-square ($m \times n$) Vandermonde matrix with a number of additional rows.
This gives an underdetermined system without a unique solution for the
coefficients $\boldsymbol{\alpha}$.  We select a solution that minimizes
$\|\boldsymbol{\alpha}\|_{2}$ for further study:  
\begin{align}
V^{T} \cdot \boldsymbol{\alpha} = \int_0^{\Delta t} \tau^{i} d\tau, \hspace{2mm} \|\boldsymbol{\alpha}\|_{2} \rightarrow \text{min}, \quad V = \begin{bmatrix} 
    1 & t_{1} &  \hdots   & t_{1}^{n-1}  \\
    1 & t_{2} & \hdots &  t_{2}^{n-1} \\
      \vdots  & \vdots  &  \ddots   &  \vdots   \\
       1 &   t_{m}  &  \hdots  & t_{m}^{n-1} \\   
        \end{bmatrix}. \label{eq:ext_ab}
\end{align}
We expect that minimization of the 2-norm of the coefficients
$\boldsymbol{\alpha}$ to be used in extrapolation will produce schemes with
improved stability characteristics.  To further motivate this choice of
$\boldsymbol{\alpha}$, we use a model ODE:
\begin{align}
\frac{dy}{dt} = -\beta y, \quad y(0) = 1. \notag
\end{align}
Referencing \eqref{eq:ab_general}, we can derive a step matrix for a third-order Adams-Bashforth integrator applied to this case such that
\begin{align}
\vec{\phi}_{\text{n+1}} = G\vec{\phi}_{\text{n}}, \notag
\end{align}
where the vector $\vec{\phi}_{n}$ is given by $\vec{\phi}_{n} = [y(t_{n})
\hspace{2mm} y(t_{n-1}) \hspace{2mm} y(t_{n-2})]^{T}$, $\vec{\phi}_{n+1}$ is
given by $\vec{\phi}_{n+1} = [y(t_{n+1}) \hspace{2mm} y(t_{n}) \hspace{2mm}
y(t_{n-1})]^{T}$, and the step matrix $G$ is given by
\begin{align}
G = \begin{bmatrix} 
    1 - \beta \alpha_{1} & -\beta \alpha_{2} & -\beta \alpha_{3}  \\
    1 & 0 & 0 \\
       0 &   1  &   0 \\   
        \end{bmatrix}, \notag
\end{align}
where $\boldsymbol{\alpha}$ is the vector of Adams-Bashforth coefficients
$\boldsymbol{\alpha} =[\alpha_{1} \hspace{2mm} \alpha_{2} \hspace{2mm}
\alpha_{3}]^{T}$.  Based on the definition of $\|G\|_{2}$, it follows that the
Adams-Bashforth integration of this case will remain stable if
$\lambda_{i}(G^{T}G) \leq 1$, where $\lambda_{i}(G^{T}G)$ are the eigenvalues
of $G^{T}G$.  The eigenvalues of $G^{T}G$ are given by
\begin{align}
\lambda_{1} = 1, \quad \lambda_{2,3} = \frac{1}{2}(C \pm \sqrt{(-C)^{2} - 4\alpha_{3}^{2}\beta^{2}}), \notag
\end{align}
where
\begin{align}
C = - 2\alpha_{1}\beta + 2 + \beta^{2}\|\boldsymbol{\alpha}\|_{2}^{2}. \notag
\end{align}
Therefore, as $\|\boldsymbol{\alpha}\|_{2} \rightarrow 0$, $\lambda_{2,3}
\rightarrow 0$ for all real values of $\beta$.  This indicates for this case
that a choice of $\boldsymbol{\alpha}$ that minimizes
$\|\boldsymbol{\alpha}\|_{2}$ is likely to produce superior stability
characteristics, as it will lead to lower values of $\|G\|_{2}$.
\revanswerb{This outcome for the scalar ODE example motivates an application of
extended-history Adams-Bashforth integrators using this coefficient choice to
systems described by the SBP-SAT discretization of the Navier-Stokes equations given in 
Section \ref{background}.  We leave a rigorous demonstration of the improved
stability of these schemes  relative to standard Adams-Bashforth schemes of the
same order on these systems for future work.}

\begin{figure}[h!]
\captionsetup[subfigure]{width=0.9\textwidth}
\centering
\begin{subfigure}[b]{0.5\textwidth}
 \centering
  \includegraphics[width=0.9\linewidth]{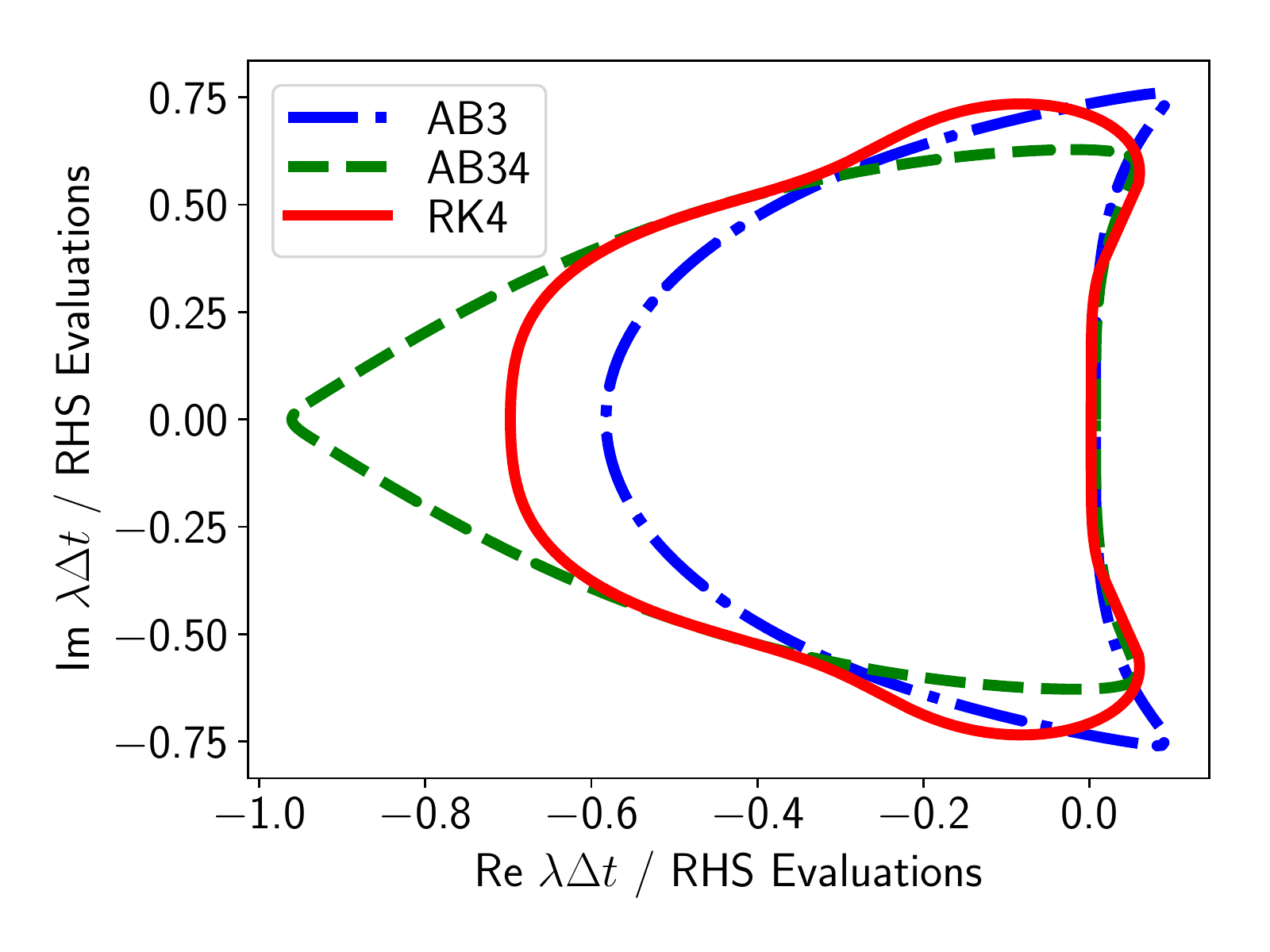}
  \caption{Approximate stability regions for AB34 integrators compared to AB3 and RK4.}
  \label{fig:stabplots34}
\end{subfigure}\begin{subfigure}[b]{0.5\textwidth}
 \centering
  \includegraphics[width=0.9\linewidth]{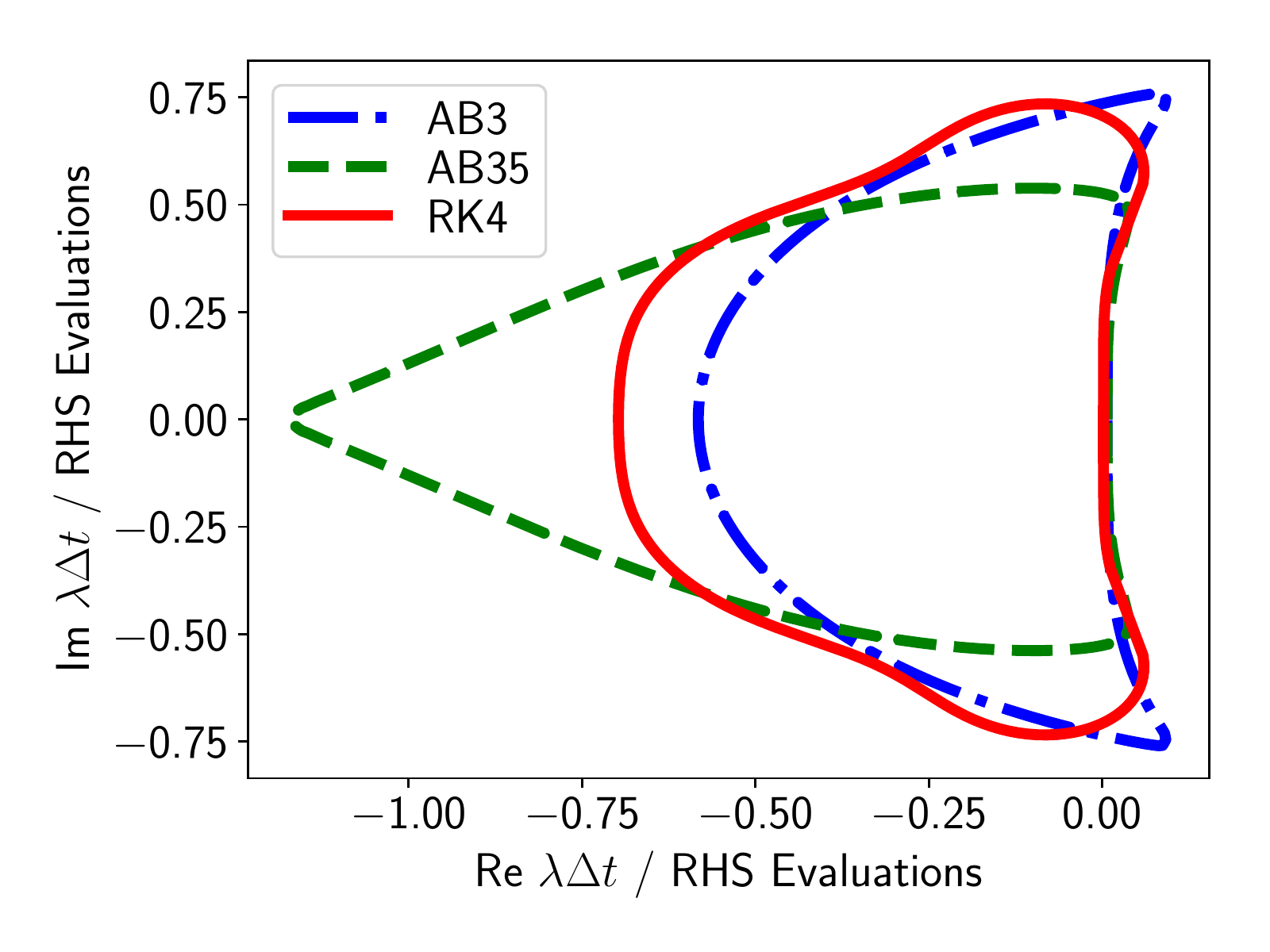}
  \caption{Approximate stability regions for AB35 integrators compared to AB3 and RK4.}
  \label{fig:stabplots35}
\end{subfigure}
\caption{Approximate stability regions obtained via the method of \ref{approx_stab} for AB3 integrators with extended history, normalized by the number of RHS evaluations.}
\label{fig:stab_plots_ext_3}
\end{figure}

\begin{figure}[h!]
\captionsetup[subfigure]{width=0.9\textwidth}
\centering
\begin{subfigure}[b]{0.5\textwidth}
 \centering
  \includegraphics[width=0.9\linewidth]{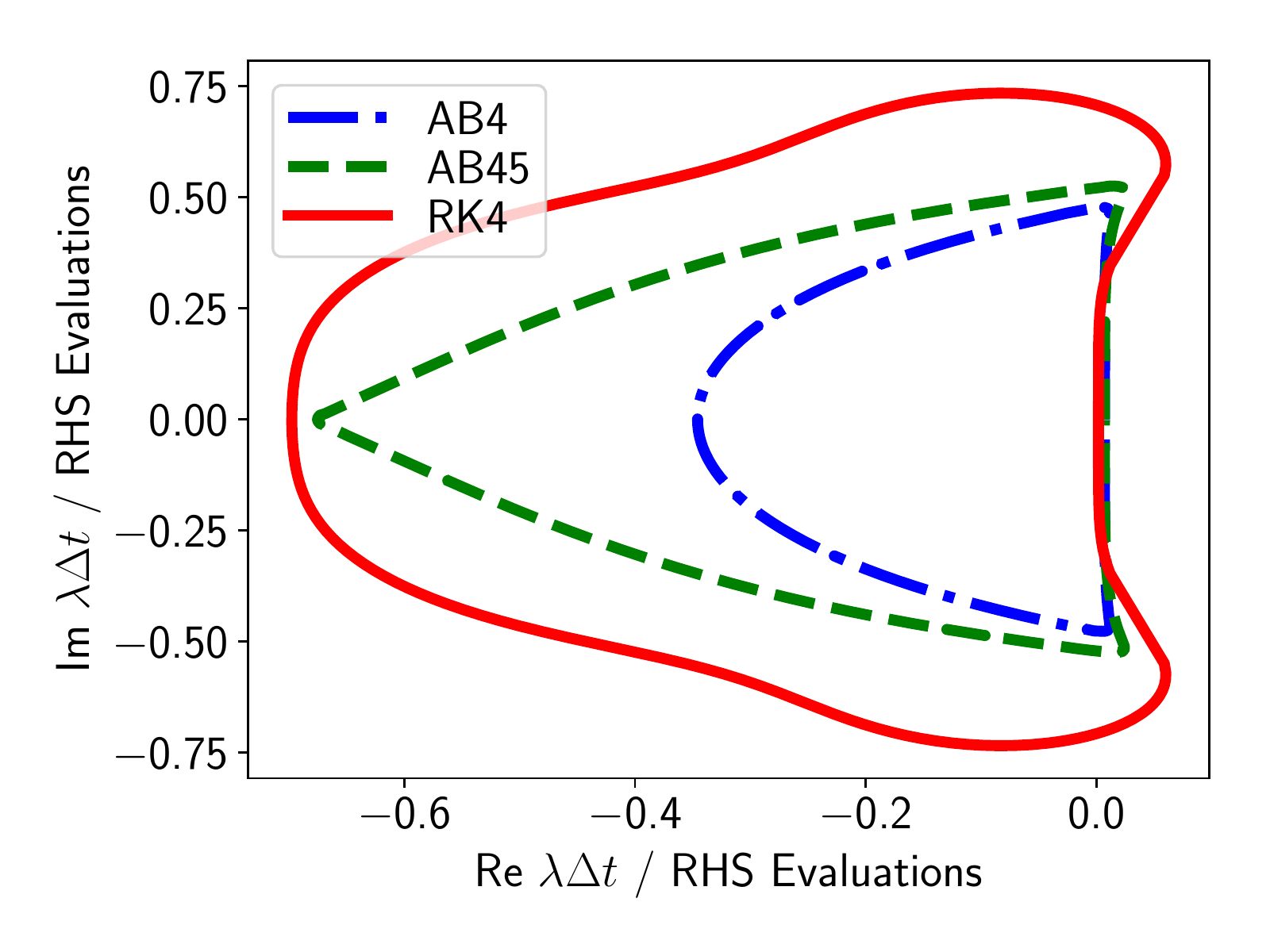}
  \caption{Approximate stability regions for AB45 integrators compared to AB4 and RK4.}
  \label{fig:stabplots45}
\end{subfigure}\begin{subfigure}[b]{0.5\textwidth}
 \centering
  \includegraphics[width=0.9\linewidth]{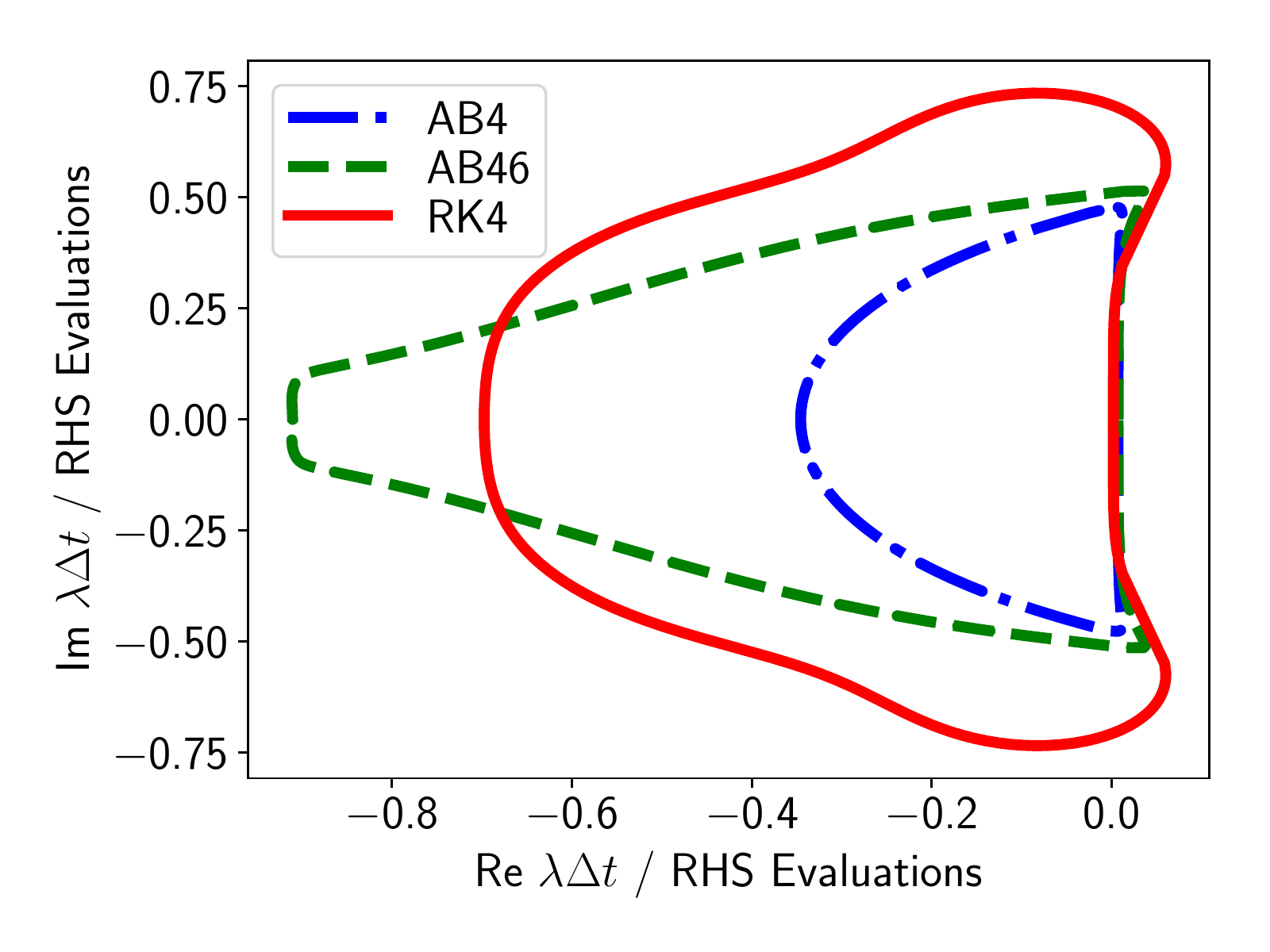}
  \caption{Approximate stability regions for AB46 integrators compared to AB4 and RK4.}
  \label{fig:stabplots46}
\end{subfigure}
\caption{Approximate stability regions obtained via the method of \ref{approx_stab} for AB4 integrators with extended history, normalized by the number of RHS evaluations.}
\label{fig:stab_plots_ext_4}
\end{figure}

This method can be extended to any order and any number of additional history
nodes.  For third and fourth-order, the method produces what we term AB34
(third-order) and AB45 (fourth-order) schemes.
Figures~\ref{fig:stab_plots_ext_3} and \ref{fig:stab_plots_ext_4} plot
stability regions (using the method we describe in \ref{approx_stab}) for the
new AB34 and AB35 schemes for comparison against a fourth-order Runge-Kutta
(RK4) integrator and standard AB3 integrators and, similarly, of the new AB45
and AB46 schemes compared to RK4 and fourth-order Adams-Bashforth (AB4)
integrators.  Motivated by the amount of computational work per timestep, we
normalize the approximate stability regions based on the number of right-hand
side evaluations per timestep a given integrator requires --- that is, while
RK4 requires four right-hand side evaluations per timestep, Adams-Bashforth
integrators only require one.

When compared with the approximate stability regions for standard
Adams-Bashforth schemes, the approximate stability regions for the
AB34 and AB45 schemes extend further along the negative real axis, while
slightly shrinking along the imaginary axis.  In the upcoming section, we will
see that the SBP-SAT discretization we use motivates the use of this scheme to
handle eigenvalues associated with inter-grid fluxes and,
in the results section, we will therefore consider both the standard AB schemes
(AB3, AB4) and the new AB34 and AB45 schemes in our analysis of convergence,
stability, and performance considerations.  While we will limit our continued
discussion of extended-history schemes to AB34 and AB45 integrators, we can
also extend the history used to produce the coefficients $\boldsymbol{\alpha}$
even further, creating AB35 and AB46 schemes, and again viewing normalized
approximate stability regions for these schemes (Figures \ref{fig:stabplots35}
and \ref{fig:stabplots46}).

\subsection{Multi-rate Adams-Bashforth Integration} \label{mrab}

We now describe a multi-rate generalization of the scheme, making use of the
algorithm introduced in~\cite{gear1984multirate}.  We consider the following
model system with ``fast'' and ``slow'' solution components:
\begin{align}
\frac{\textrm{d}}{\textrm{d}t}\left( \begin{array}{c} f(t) \\ s(t) \end{array} \right) = \left( \begin{array}{c} a_{f}(f,s) \\ a_{s}(f,s) \end{array} \right). \label{eq:mrab}
\end{align}
With this in mind, we can set a slow (larger) time step $H$ for $a_{s}$ such
that we maintain stability in the integration of the slow component.  We also
set a fast time step $h$ for $a_{f}$ such that $H$ is an integer multiple of
$h$, and define the ratio between the two, $\text{SR} = H/h$, assumed to be an
integer, as the \emph{step ratio} of the MRAB scheme. While the results presented here
make use of only two separate state components, each with its own right-hand
side function and independent rate, the approach is readily extensible to any
number of rates.

In the two-component overset formulation with which we are concerned, we define
the fast and slow components of our Navier-Stokes solution as the conserved
variables $Q_{i} = [\rho_{i}, (\rho \vec{u})_{i}, (\rho E)_{i}]^T$ on each
grid, that is (using a two-grid case as an example): $f = Q_{1}$, $s = Q_{2}$,
where the subscripts of the vectors $Q$ indicate global grid number.  
We assume Grid 1 to be the grid with the fast-moving component of the solution, 
be it due to physical behavior or finer mesh spacing.  Each
right-hand side function $a_{f}$ and $a_{s}$ is a function of both the slow and
fast states $s$ and $f$ --- this coupling between the right-hand side functions
is, in the case of our application of this theory to overset meshes, realized
by the SAT penalty interpolation discussed in Section~\ref{interp}.

Within this two-component scheme, a few design choices are available:
\begin{itemize}
\item The order in which we evaluate and advance the solution components.
	Namely, two primary options are advancing the fast-evolving solution
	component through all of its micro-timesteps $h$ and waiting to perform
	the single macro-timestep $H$ required for the slow component until the
	end (a ``fastest-first'' scheme, per the nomenclature of
	\cite{gear1984multirate}), or pursuing an algorithm in which the slow
	component is instead advanced first.
\item For slowest-first evaluation schemes, the choice of whether or not to
	re-extrapolate the slow state after additional state and right-hand
	side information is gathered at the micro-timestep level.
\end{itemize}
Empirical observations on the effects of these choices are made
in \citep{klockner2010high}.  It is useful to step through a brief example
of a multi-rate Adams-Bashforth integrator, using a system with a
fast component requiring twice as many timesteps as the slow component
to remain well-resolved ($\text{SR}=2$).  We lay out the steps of a third-order
fastest-first MRAB scheme with no re-extrapolation, assuming that $a_{s}$
evolves at the slow rate (macro-timestep $H = 2h$) and $a_{f}$ evolves at the
fast rate (micro-timestep $h$).  $\hat{a}$ denotes extrapolants of the
right-hand side functions as polynomial functions of both the set of time history points $\vec{t}$ and the set
of history values $\vec{a}_{\text{hist}}$: $\hat{a} = P(t, \vec{t},
\vec{a}_{\text{hist}})$, where
\begin{align}
P(t, \vec{t}, \vec{a}_{\text{hist}}) = \sum_{k=1}^{n} \Big(\prod_{\substack{ 0 \leq j \leq n \\ j \neq k}} \frac{t - t_{j}}{t_{k} - t_{j}} \Big) a_{\text{hist},k} \notag
\end{align}
These polynomials approximating the evolution of
$a_{f}$ and $a_{s}$ in time are what we will integrate to march $a_{f}$ and
$a_{s}$, and they will be updated to replace older history values with new
right-hand side evaluations during the course of integration through a
macro-timestep $H$.  We assume availability of right-hand side histories to
start the AB method.
\begin{itemize}[leftmargin=0.5in]
\item[Step 1:] Form the polynomial extrapolants we will integrate, per the AB
	methods described in Section~\ref{ab}:	
\begin{align}
	\hat{a}_{f,1}(t) = P\Big(t, [t_{i-2}, t_{i-1}, t_{i}], [a_{f}(f(t_{i-2}), s(t_{i-2})), a_{f}(f(t_{i-1}), s(t_{i-1})), a_{f}(f(t_{i}), s(t_{i}))]\Big) \notag \\
	\hat{a}_{s,1}(t) = P\Big(t, [t_{i-4}, t_{i-2}, t_{i}], [a_{s}(f(t_{i-4}), s(t_{i-4})), a_{s}(f(t_{i-2}), s(t_{i-2})), a_{s}(f(t_{i}), s(t_{i}))]\Big) \notag
\end{align}	
	The right-hand side history
	values of $a_{f}$ (used to form $\hat{a}_{f,1}$) have been obtained at
	time points $t_{i-2} = t - 2h_{i}$, $t_{i-1} = t - h_{i}$, and current time
	$t_{i} = t$, whereas the right-hand side history values of $a_{s}$
	(used to form $\hat{a}_{s,1}$) have been obtained at time points $t_{i-4}
	= t - 2H_{i}$, $t_{i-2} = t - H_{i}$, and $t_{i} = t$.  The macro-timestep $H_{i}$
	can change on a per-macrostep basis, with the micro-timestep $h_{i}$ defined such that
	$h_{i} = H_{i}/\text{SR}$.
\item[Step 2:] March both $f$ and $s$ to time $t_{i+1}$ by integrating the
	polynomial extrapolants $\hat{a}_{s,1}(t)$ and $\hat{a}_{f,1}(t)$ formed in
	Step~1:
\begin{align}
f(t_{i+1}) = f(t_{i}) + \int_{t_{i}}^{t_{i+1}} \hat{a}_{f,1}(\tau) d\tau, \notag \\
s(t_{i+1}) = s(t_{i}) + \int_{t_{i}}^{t_{i+1}} \hat{a}_{s,1}(\tau) d\tau. \notag
\end{align}
This results in a set of intermediate values $f(t_{i+1})$ and $s(t_{i+1})$.
\item[Step 3:] Evaluate the fast right-hand side $a_{f}(f(t_{i+1}), s(t_{i+1}))$.
\item[Step 4:] Update the set of right-hand side history values for $a_{f}$ to include these new values, and construct a new extrapolant $\hat{a}_{f,2}$:
\begin{align}
\hat{a}_{f,2}(t) = P\left(t, [t_{i-1}, t_{i}, t_{i+1}], [a_{f}(f(t_{i-1}), s(t_{i-1})), a_{f}(f(t_{i}), s(t_{i})), a_{f}(f(t_{i+1}), s(t_{i+1}))]\right). \notag
\end{align}
\item[Step 5:] March $s$ to time $t_{i+2}$ by integrating the extrapolant formed in Step~1:
\begin{align}
s(t_{i+2}) = s(t_{i}) + \int_{t_{i}}^{t_{i+2}} \hat{a}_{s,1}(\tau) d\tau. \notag
\end{align}
\item[Step 6:] March $f$ to time $t_{i+2}$ by integrating the extrapolant formed in Step~3:
\begin{align}
f(t_{i+2}) =  f(t_{i+1}) + \int_{t_{i+1}}^{t_{i+2}} \hat{a}_{f,2}(\tau) d\tau. \notag
\end{align}
\item[Step 7:] Go to Step 1.
\end{itemize}
The scheme evaluates the fast-evolving right-hand side $a_{f}$ twice per
macro-timestep $H$, whereas the slowly-evolving right-hand side $a_{s}$ is only
evaluated once.  For the results shown later, this is the scheme we will use,
albeit generalized to different step ratios $\text{SR} = H/h$.  

\section{Time Integration in SBP-SAT Discretizations} \label{ns_stability}

One of the critical questions to be answered is at what timestep sizes $h$ and
$H$ the integrators developed in Section~\ref{integration} remain stable.  This
section establishes the procedures and tools for characterizing the stability
of the integrators developed in Section \ref{integration} in the context of the
SBP-SAT discretization used.

\subsection{Timestep Calculation} \label{ts}
For the discretization of the Navier-Stokes equations described in
Section~\ref{background}, we will incorporate the metrics and
multidimensionality into calculation of the timestep used by a given
integrator.  The maximum stable timestep for a given Navier-Stokes simulation
using the SBP-SAT discretization is approximated via the following steps.
First, we calculate an inviscid timestep:
\begin{align}
\Delta t_{\text{inv}} = \frac{J^{-1}}{\sum\limits_{i=1}^{N_{D}}|\vec{U} \cdot \frac{1}{J}\frac{\partial \xi_{i}}{\partial \vec{x}}| + c \sqrt{\hat{\xi_{k}} \cdot \hat{\xi_{k}}}} , \label{eq:ts_inv}
\end{align}
and a viscous timestep:
\begin{align}
\Delta t_{\text{visc}} = \frac{J^{-1}}{\sum\limits_{i=1}^{N_{D}}|\vec{U} \cdot \frac{1}{J}\frac{\partial \xi_{i}}{\partial \vec{x}}| + c \sqrt{\hat{\xi_{k}} \cdot \hat{\xi_{k}}} + 2 \nu^{*} J \sum\limits_{i=1}^{N_{D}}\sqrt{\frac{1}{J}\frac{\partial \xi_{i}}{\partial \vec{x}} \cdot \frac{1}{J}\frac{\partial \xi_{i}}{\partial \vec{x}}})^{2}}, \label{eq:ts_visc}
\end{align}
where the vector $\vec{U}$ contains the Cartesian velocities, $c$ is the speed
of sound, $\nu^{*} = \max{(\mu, k/c_{v})}$, $J$ is the Jacobian, $N_{D}$ is the
number of spatial dimensions, and the $\xi$ terms are the metric terms
described in \cite{maccormack2014numerical}.  To determine the timestep, we
calculate $\Delta t_{\text{inv}}$ and $\Delta t_{\text{visc}}$ at each mesh
point, then take
\begin{align}
\Delta t = \min{(\Delta t_{\text{inv}}, \Delta t_{\text{visc}})}, \label{eq:ts_overall}
\end{align}
and take an additional minimum over all mesh points.  This timestep calculation is described in more detail in \cite{maccormack2014numerical}.

\subsection{SAT Interpolation Stability} \label{sat_stab}
We now undertake a brief parameter study to show that one of the stiffest
right-hand side components limiting the timestep when solving the Navier-Stokes
equations with the overset discretization involves the SAT penalty-based
interpolation (discussed in Section~\ref{interp}).  Specifically, we will show
that penalty terms applied at fringe points with more than one nonzero
characteristic direction (i.e. corner points) produce eigenvalues with large
negative real part.  To illustrate this phenomenon, we can linearize
a computationally small Navier-Stokes problem (using the procedure for
approximation of the Jacobian described in \cite{osusky2010parallel} to form
the approximate linear operator of the system) employing SAT interpolation on
overset meshes and analyze the spectrum of the resulting global operator.  

We express the time evolution of a fringe corner point as 
\begin{align}
\frac{\partial \mathbf{q}}{\partial t} = \text{RHS}_{\text{NS}}(\mathbf{q}) +  \Omega \cdot \text{RHS}_{\text{INT}}(\mathbf{q},\hat{\mathbf{q}}), \label{eq:penalty_interp}
\end{align}
with $\text{RHS}_{\text{NS}}(\mathbf{q})$ denoting the Navier-Stokes right-hand
side, and $\text{RHS}_{\text{INT}}(\mathbf{q},\hat{\mathbf{q}})$ denoting the
penalty terms applied by the SAT interpolation procedure described in Section
\ref{interp} (the latter two terms in \eqref{eq:int}).  We introduce a
parameter $\Omega \in [0,1]$ to change the magnitude of terms associated with
the interpolation penalty applied at the corner fringe points, and determine
the right-hand side modes associated with these terms.

\begin{figure}[h!]
\captionsetup[subfigure]{width=0.9\textwidth}
\centering
\begin{subfigure}[b]{0.5\textwidth}
 \centering
  \includegraphics[width=0.9\linewidth]{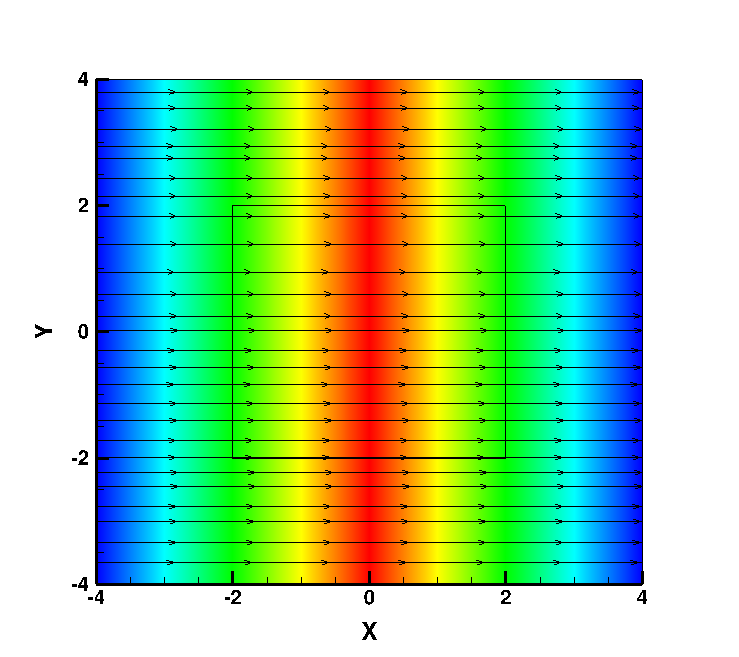}
  \caption{Penalty parameter study - 1D Euler test case initial conditions.}
  \label{fig:bib_case}
\end{subfigure}\begin{subfigure}[b]{0.5\textwidth}
 \centering
  \includegraphics[width=0.9\linewidth]{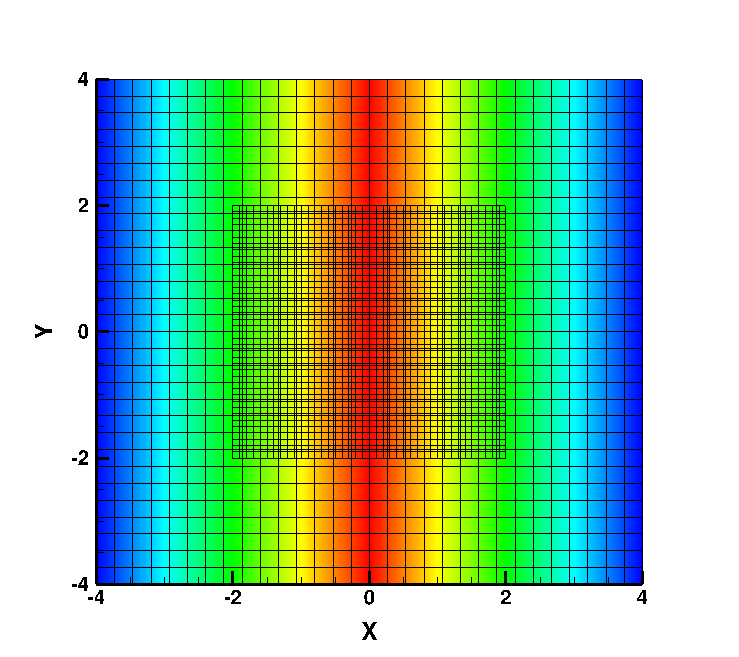}
  \caption{Penalty parameter study - 1D Euler test case meshes.}
  \label{fig:bibmesh}
\end{subfigure}
\caption{Description of the (effectively) 1D Euler test case for evaluating corner interpolation.}
\label{fig:bib_grids}
\end{figure}

As for the case to be analyzed, we select an inviscid application of overset
grids to the Euler equations with periodic boundary conditions.  The initial
condition is an $x$-momentum gradient of the form
\begin{align}
\rho U = 1.25 - 0.01|x|, \notag
\end{align}
with the $y$-momentum set to 0 everywhere. The initial condition and grid
layout for this case are shown in Figure~\ref{fig:bib_grids}.  Grid 1 ($31
\times $31 = 961 points) spans $x = [-4,4]$, $y = [-4,4]$, while Grid 2 ($41
\times 41$ = 1681 points) spans $x = [-2,2]$, $y = [-2,2]$.

\begin{figure}[h!]
\captionsetup[subfigure]{width=0.9\textwidth}
\centering
\begin{subfigure}[b]{0.5\textwidth}
 \centering
  \includegraphics[width=0.94\linewidth]{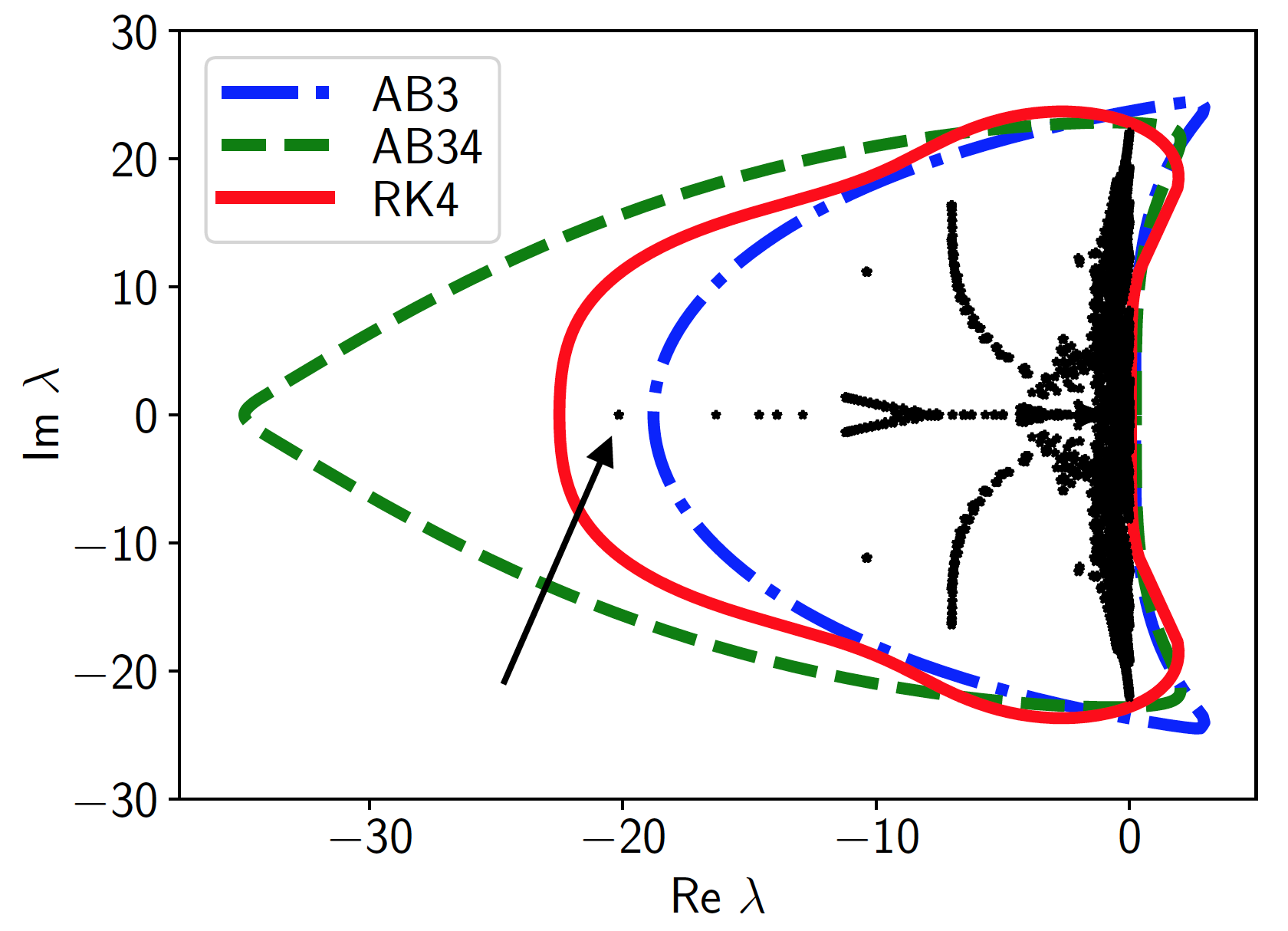}
  \caption{Spectrum for linearized operator with $\Omega = 1$, including AB34 scheme.}
  \label{fig:fullcorner_ab34}
\end{subfigure}\begin{subfigure}[b]{0.5\textwidth}
 \centering
  \includegraphics[width=0.95\linewidth]{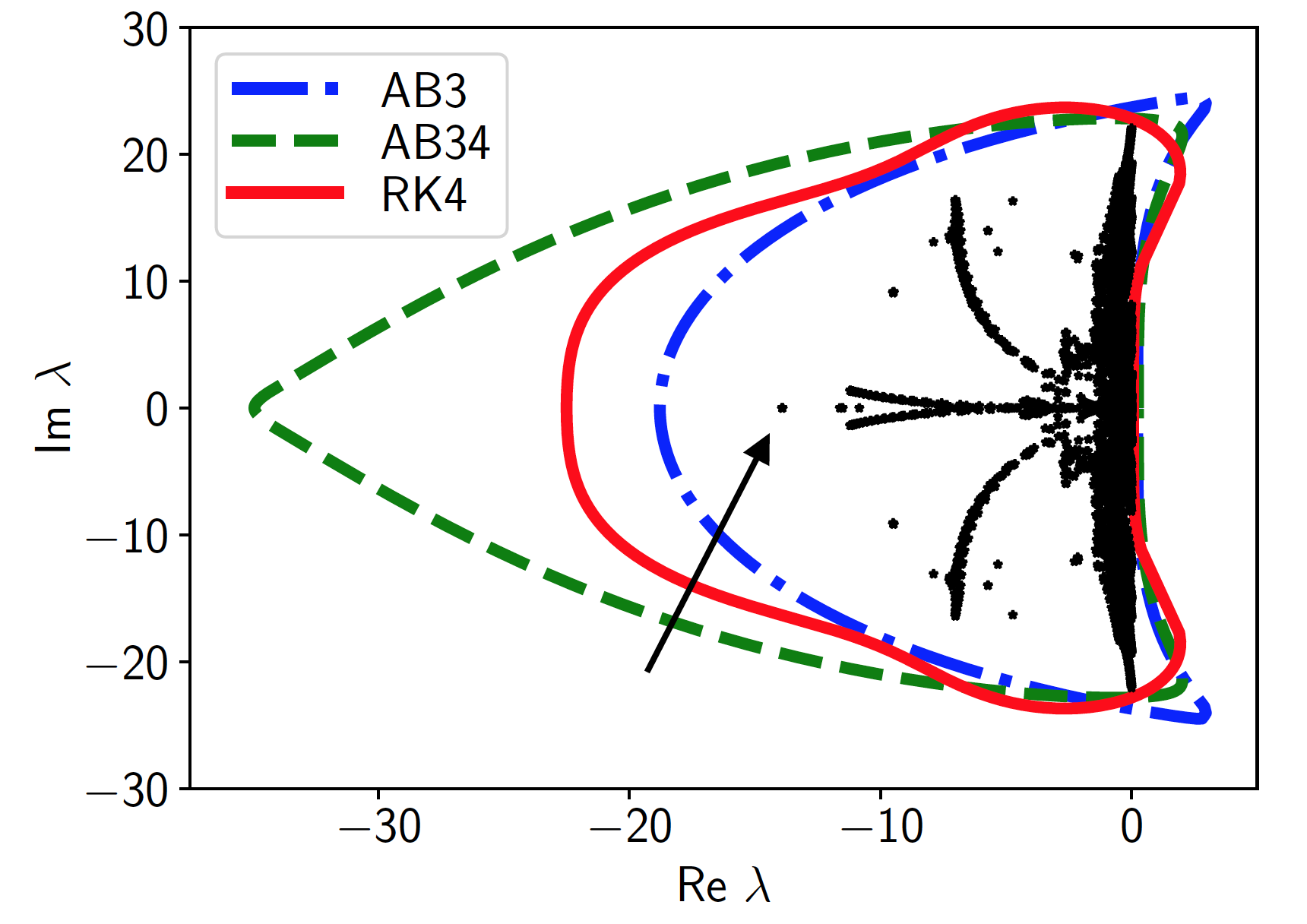}
  \caption{Spectrum for linearized operator with $\Omega = 0.5$, including AB34 scheme.}
  \label{fig:halfcorner_ab34}
\end{subfigure}
\caption{Spectra of box-in-box case for the linearized operator with $\Omega = 1, 0.5$, including AB34 scheme.  Stability regions are plotted using the method described in \ref{approx_stab}.}
\label{fig:corner_treatments_ab34}
\end{figure}

We form the approximate linearized operators for this case for two values of
the corner penalty weight parameter $\Omega$: $\Omega = 1$ and $\Omega = 0.5$.
By finding the eigenvalues of the two operators and comparing them, we can
empirically identify the eigenvalues associated with the corners.  Based on the
highest-magnitude corner-associated eigenvalue that can be determined from
comparing these spectra, it is clear that the approximate stability region of
the RK4 scheme contains the corner-associated eigenvalue, which extends
outwards along the negative real axis as we increase the penalty weight
$\Omega$ from 0.5 to 1, whereas the approximate stability region for the
third-order AB scheme does not, where in each case the time step was calculated
according to \ref{eq:ts_overall}, and in which the scaling is controlled by the
behavior along the imaginary axis in the situation considered.  This motivated
the development of the extended-history schemes discussed in
Section~\ref{ab_new}, which are shown to improve the stability of AB schemes
along the negative real axis.  Figure~\ref{fig:corner_treatments_ab34} shows
that the approximate stability region for the third-order extended-history AB
scheme indeed contains the corner-associated eigenvalue.

\subsection{Procedure for Determining Discrete Stability} \label{sm}
We now develop a procedure for characterizing the discrete temporal stability
of a given integrator within the SBP-SAT discretization using the same step
matrix formulation as in the analysis of Section~\ref{ab_new}.  We can
express integration forward one timestep as a matrix operation on the
integrator state vector $\vec{\phi}_{\text{n}}$:
\begin{align}
\vec{\phi}_{\text{n+1}} \approx G_{\epsilon}\vec{\phi}_{\text{n}}, \label{eq:sm}
\end{align}
where $G_{\epsilon}$ denotes the ``step matrix'' for a given integrator.  These
step matrices are functions of the timestep $H$ used in the simulation, and are
formed using a linearization procedure that approximates the full Jacobian
matrix of the system using finite differences.  In doing so, analyzing the
eigenvalues of this step matrix is analogous to analyzing the amplification
factors of various modes in the global error of the scheme.

To construct the matrix $G_{\epsilon}$, we first use an initial state as
input to an integrator with timestep $H$ to march one step forward in
time and form the base result about which we linearize.  Next, we perturb a single element of the
initial state vector using a small value $\epsilon$.  Using this perturbed
state as an input to the integrator returns a state that --- when compared to
the base result --- demonstrates the effect of the perturbation of that single
element on the entire domain.  If we define $\vec{e}_{j}$ as a column vector
with all zeroes except for a 1 in the $j$th row, we can write
\begin{align}
	\vec{\phi}_{\text{n+1, base}} \approx G_{\epsilon}\vec{\phi}_{\text{n, init}} \notag \\
\vec{\phi}_{\text{n, mod}} = \vec{\phi}_{\text{n, init}} + \epsilon \vec{e}_{j} \notag \\
\vec{\phi}_{\text{n+1, mod}} \approx G_{\epsilon}\vec{\phi}_{\text{n, mod}}. \notag
\end{align}
Thus we calculate the $j$th column of the step matrix by perturbing the $j$th
initial state vector element, giving this modified initial state to the
integrator to obtain $\vec{\phi}_{\text{n+1, mod}}$, and taking
\begin{align}
G_{\epsilon}(:,j) \approx \frac{\vec{\phi}_{\text{n+1, mod}} - \vec{\phi}_{\text{n+1, base}}}{\epsilon}. \label{eq:sm_form}
\end{align}
Note also that in the case of Adams-Bashforth integration, the stability vector
element being perturbed may also be a right-hand side history element.  For our
perturbation $\epsilon$, we use a value of $10^{-7}$.  This procedure is
similar to the Navier-Stokes linearization for SBP-SAT discretizations
described in \cite{osusky2010parallel}.  

Once the matrix $G_{\epsilon}$ is formed, we consider the system discretely
stable if the spectral radius $\rho(G_{\epsilon})$ fulfills the property
$|\rho(G_{\epsilon})| \leq 1$.  Therefore, creating step matrices for an
integrator at a number of different values of $\Delta t$ will allow us to
characterize that integrator's stability.  We use SLEPC
\cite{Hernandez:2005:SSF} for the computation of the eigenvalues, using the
default Krylov-Schur solver with a tolerance of $10^{-12}$.  We use this
procedure to determine the maximum stable timesteps of our integrators reported
in Section \ref{stab_results}.

\subsection{Validation of the approximate stability calculation via RK4} \label{rk_stab}

We validate the procedure of the previous section by using it to determine the
maximum stable timestep for an explicit fourth-order Runge-Kutta (RK4)
integrator applied to 1D advection with a wave speed of 1.  With a spatial
domain $x = [0,1.01666]$ with 61 mesh points, the step matrix analysis
described in Section~\ref{sm} predicts a maximum stable timestep for this case
between $\Delta t = 0.03$ and $\Delta t = 0.035$. Using the method described in
\ref{approx_stab} to approximate a given integrator's stability region in the
complex plane, we see that the imaginary-axis bound for an RK4 integrator is
about $\lambda \Delta t_{\text{max}} = 2.78 i$.  When using a fourth-order SBP
operator as described in \citep{strand1994summation, carpenter1993time,
mattsson2004stable} for spatial discretization of the semidiscrete problem, and
with periodic boundaries, the eigenvalues can be found analytically, and all
lie on the imaginary axis, with a maximum eigenvalue of $\lambda_{\text{max}}
\Delta x= 1.3722 i$ (a result given by Lele~\cite{lele1992compact}).  We can
thus estimate the maximum stable timestep for this problem, noting a grid
spacing of $\Delta x = 0.01666$, to be $\Delta t_{\text{max}} = \frac{\lambda
\Delta t_{\text{max}}}{\lambda_{\text{max}} \Delta x} \Delta x = 0.03377$, a
result which is consistent with that of our step matrix procedure.

\section{Numerical Results} \label{results}

\subsection{Guiding Experiment: Viscous Flow Past a Cylinder}  \label{case}

We consider two-dimensional viscous flow over a cylinder with diameter $D =
0.6$ in the spatial domain $[-4,4]^{2}$.  The cylinder center is
located at $x = -1.2$, $y = 0$.  The free-stream Mach number is 0.2, the
Reynolds number is 200 and the Prandtl number is 0.72.  We start our
simulations from a steady state solution obtained after 100 non-dimensional
time units ($\text{``NDTU''} = tU_{\infty}/D$, using the nondimensionalization
of \cite{park1998numerical}) using an RK4 integration scheme.  Note that the
timesteps we report are also non-dimensional.  The physical problem is modeled
using two overset meshes: a coarser, base Cartesian grid ($61 \times 61$), and
a finer teardrop-shaped curvilinear grid ($121 \times 40$) surrounding the
stationary cylinder.  

\begin{figure}[h!]
\centering
\begin{subfigure}[b]{0.5\textwidth}
 \centering
  \includegraphics[width=0.9\linewidth]{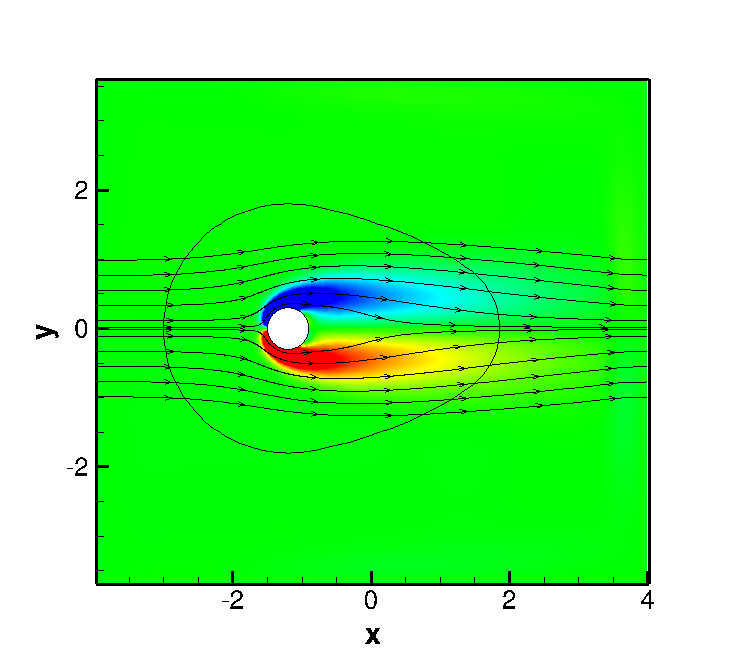}
  \caption{Physical case: flow past a cylinder.}
  \label{fig:cyl1}
\end{subfigure}\begin{subfigure}[b]{0.5\textwidth}
 \centering
  \includegraphics[width=0.9\linewidth]{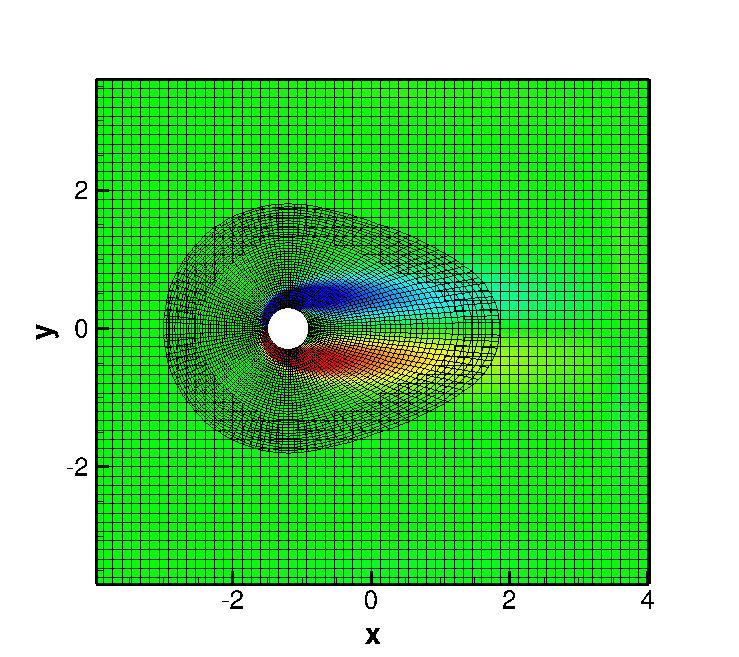}
  \caption{Overset grids: flow past a cylinder.}
  \label{compmodel}
\end{subfigure}
\caption{Physical and computational description of the flow past cylinder test case used to obtain numerical results.}
\label{fig:cylinder}
\end{figure}

It is important to characterize the disparity in allowable timestep sizes such
that we can estimate the maximum allowable step ratio $\text{SR}$ we are able
to use in our integrators.  Calculated using the method of \eqref{eq:ts_inv},
\eqref{eq:ts_visc}, and \eqref{eq:ts_overall}, we find the ratio between the
maximum timesteps on each grid to be about 12.

The boundary conditions employed in our simulation are as follows:  The inner
fine grid (Grid 2) is periodic in the azimuthal direction.  We model the
cylinder surface as an SAT isothermal wall per \cite{svard2008stable} (Equation 19).  On Grid
1 (Cartesian base grid), all outer boundaries are modeled as SAT far-field
following the procedure of \cite{svard2007stable} (Equation 22).  In addition, sponge
boundaries with a cell depth of 6 are used as described in
\cite{bodony2006analysis} (Equation 69).

\subsection{Multi-rate Adams-Bashforth Stability} \label{stab_results}

We will now determine the stability characteristics of a given MRAB integrator
by defining two ratios: $r_{\text{RK4}} = \Delta t/\Delta t_{\text{RK4}}$ and
$r_{\text{SRAB}} = \Delta t/\Delta t_{\text{SRAB}}$, where $\Delta t$ is the
maximum stable timestep of the integrator in question, and $\Delta
t_{\text{SRAB}}$ and $\Delta t_{\text{RK4}}$ are the maximum stable timesteps
for single-rate Adams-Bashforth (SRAB) and RK4 integrators, respectively,
determined using the method of Section \ref{sm}.  Recall that the step ratio of
a given MRAB integrator is given as $\text{SR}=H/h$, and is an integer value.
An MRAB integrator with a step ratio of 1 is equivalent to an SRAB integrator.
$\text{SR}$ is a user-specified input.

For the experiments to determine stability, we set the value of the
macro-timestep $H = \Delta t$, build the step matrix of \eqref{eq:sm} (starting
from a steady-state solution), and check for eigenvalues outside the unit disc,
repeating the procedure until we find the maximum critical value of $\Delta t$
(to the nearest $10^{-4}$) for an MRAB integrator with a given step ratio.  The
$r_{\text{RK4}}$ values are a measure of the maximum stable timesteps of these
integrators relative to an RK4 integrator.  The $r_{\text{SRAB}}$ values
reported for an MRAB integrator with a given order and history length use the
maximum stable timestep for an SRAB integrator of the same order and history
length.  The percentages reported document percent efficiency of a given
integrator, $\% = r_{\text{SRAB}}/\text{SR}$: that is, if an MRAB integrator
with a step ratio of $\text{SR}$ obtains an $r_{\text{SRAB}}$ value of SR, it
is deemed 100\% efficient.

\begin{table}[h] \begin{center}
  \begin{tabular}{cccccccccc}
       \toprule
       \multicolumn{2}{c}{} & \multicolumn{4}{c}{Standard (AB3)} & \multicolumn{4}{c}{Extended-History (AB34)} \\
       \cmidrule(lr){1-2}  \cmidrule(lr){3-6}  \cmidrule(lr){7-10}
       Int. &                     SR &                      $\Delta t$ &             $r_{\text{RK4}}$      & $r_{\text{SRAB}}$ &  \%   &             $\Delta t$ &             $r_{\text{RK4}}$      & $r_{\text{SRAB}}$ & \%  \\
       \midrule
     RK4 &                                  &                            0.0217 &                           1.000 &                              5.13 &                &            0.0217 &                  1.000 &             3.86  &           \\
     SRAB &                           1 &                           0.0042 &                            0.195 &                            1.00 &      100 &            0.0056 &                  0.259 &             1.00 &         100 \\
     MRAB &                           2 &                          0.0084 &                            0.388 &                            1.99 &       99.5 &             0.0113 &                  0.518 &             2.00 &          100 \\
     MRAB &                           3 &                          0.0127 &                            0.583 &                            2.99 &       99.7 &           0.0169 &                  0.775 &             2.99 &             99.7 \\
     MRAB &                           4 &                          0.0169 &                            0.777 &                            3.99 &       99.8 &           0.0225 &                 1.035 &             3.99  &             99.8 \\
     MRAB &                           5 &                          0.0171 &                            0.788 &                            4.04 &        80.8 &          0.0281 &                 1.294 &             4.99  &            99.8  \\
     MRAB &                           6 &                          0.0171 &                           0.788 &                            4.04  &         67.3 &         0.0282 &                 1.296 &             5.00  &            83.3 \\
     MRAB &                           7 &                          0.0171 &                            0.788 &                            4.04  &        57.7 &          0.0282 &                 1.296 &             5.00 &           71.4  \\
     MRAB &                           8 &                          0.0171 &                            0.788 &                            4.04 &         50.5 &          0.0282 &                 1.296 &             5.00 &           62.5  \\
     MRAB &                           9 &                          0.0171 &                            0.788 &                            4.04 &         44.9 &          0.0282 &                 1.296 &             5.00 &          55.6  \\
     MRAB &                          10 &                        0.0171 &                            0.788 &                           4.04 &          40.4 &           0.0282 &                 1.296 &             5.00 &          50.0  \\ 
     MRAB &                          20 &                       0.0171 &                            0.788 &                           4.04 &            20.2 &          0.0282 &                 1.296 &             5.00 &          25.0 \\ 
     \bottomrule                  
 \end{tabular}
 \end{center}
  \caption{Stability results for third order AB integrators applied to flow past a cylinder (speed ratio of about 12).  $r_{\text{RK4}} = \Delta t/\Delta t_{\text{RK4}}$, $r_{\text{SRAB}} = \Delta t/\Delta t_{\text{SRAB}}$, and $\% = r_{\text{SRAB}}/\text{SR}$.}  \label{t:thirdstab}
\end{table}

\begin{table}[h] \begin{center}
  \begin{tabular}{cccccccccc}
       \toprule
       \multicolumn{2}{c}{} & \multicolumn{4}{c}{Standard (AB4)} & \multicolumn{4}{c}{Extended-History (AB45)} \\
       \cmidrule(lr){1-2}  \cmidrule(lr){3-6}  \cmidrule(lr){7-10}
       Int. &                     SR &                      $\Delta t$ &             $r_{\text{RK4}}$      & $r_{\text{SRAB}}$ &  \% &             $\Delta t$ &             $r_{\text{RK4}}$      & $r_{\text{SRAB}}$ & \% \\
       \midrule
     RK4  &                                 &                         0.0217 &                            1.000 &                       9.42 &             &                   0.0217 &                 1.000 &                            5.63 &                    \\
     SRAB &                             1 &                        0.0023 &                            0.106 &                       1.00 &             100 &        0.0039 &                 0.178 &                            1.00 &          100  \\
     MRAB &                            2 &                        0.0046 &                            0.212 &                       2.00 &              100 &       0.0078 &                 0.356 &                            2.00 &           100 \\
     MRAB &                            3 &                        0.0070 &                            0.321 &                       3.02 &               100 &      0.0116 &                 0.533 &                            3.00 &            100 \\ 
     MRAB &                            4 &                        0.0093 &                            0.427 &                       4.02 &              100 &       0.0155 &                 0.711 &                            4.00 &            100 \\ 
     MRAB &                            5 &                        0.0093 &                            0.430 &                       4.05 &               81.0 &        0.0180 &                 0.827 &                            4.65 &            93.0 \\
     MRAB &                            6 &                        0.0093 &                            0.430 &                       4.05  &              67.5 &        0.0180 &                 0.827 &                           4.65 &             77.5 \\ 
     MRAB &                            7 &                        0.0093 &                            0.430 &                       4.05 &               57.9 &        0.0180 &                 0.827 &                            4.65 &            66.4 \\
     MRAB &                            8 &                        0.0093 &                            0.430 &                       4.05 &               50.6 &        0.0180 &                 0.827 &                            4.65 &            58.1 \\
     MRAB &                            9 &                        0.0093 &                            0.430 &                       4.05 &               45.0 &        0.0180 &                 0.827 &                            4.65 &            51.7  \\
     MRAB &                          10 &                        0.0093 &                            0.430 &                       4.05 &               40.5 &        0.0180 &                 0.827 &                            4.65 &            46.5 \\ 
     MRAB &                          20 &                        0.0093 &                            0.430 &                       4.05 &               20.3 &        0.0180 &                 0.827 &                            4.65 &            23.3 \\
     \bottomrule                 
 \end{tabular}
 \end{center}
  \caption{Stability results for fourth-order AB integrators applied to flow past a cylinder (speed ratio of about 12).  $r_{\text{RK4}} = \Delta t/\Delta t_{\text{RK4}}$, $r_{\text{SRAB}} = \Delta t/\Delta t_{\text{SRAB}}$, and $\% = r_{\text{SRAB}}/\text{SR}$.}  \label{t:fourthstab}
\end{table}

The results of Tables~\ref{t:thirdstab} and~\ref{t:fourthstab} show that the
extended-history schemes have higher maximum stable timesteps than their
standard counterparts of same order.  All schemes attain near-perfect
efficiency ($r_{\text{SRAB}} \approx \text{SR}$) for $\text{SR}$-values up to
about a third of the speed ratio of 12 for this case ($\text{SR} = 4$), with
the extended-history schemes maintaining this efficiency for the $\text{SR}=5$
cases as well, a result due to the improved stability of these schemes along
the negative real axis.  Above this step ratio, we expect our performance to
degrade, given that we are simply performing extra right-hand side evaluations
(more micro-timesteps on the fast grid) with no commensurate macro-timestep
gain.   The increase from third to fourth order necessitates the use of lower
timesteps for all step ratios, but the same qualitative ceiling as in the
third-order case is observed.  The third-order MRAB scheme with extended
history is the only scheme able to attain $r_{\text{RK4}} \geq 1$ for any step
ratio, indicating maximum stable timesteps larger than that of the RK4
integrator.

\subsection{Multi-Rate Adams-Bashforth Accuracy and Convergence}\label{accuracy}

We now examine the accuracy of the MRAB schemes and confirm that they are
convergent.  In presenting accuracy results for our schemes, we will
examine integrators for values of $\text{SR}$ ranging from 2 to 6.  As we have
seen in the previous section, above this step ratio, the maximum stable
timestep remains unchanged for all AB integrators for our speed ratio of about
12.  For each integrator, we calculate an estimated order of convergence (EOC)
using 3 data points consisting of the macro-timestep used and the maximum error
$\text{err} = \|\rho - \rho_{\text{true}}\|_{\infty}$ obtained in the density
after 2.5 NDTU\@.  The initial condition for these 2.5-NDTU integrations is a
uniform subsonic flow with Mach number $M = 0.2$ in the $x$-direction (before
the boundary layer has developed on the cylinder), and the true solution
$\rho_{\text{true}}$ to which the solution obtained with our integrators is
compared to form the error is the solution obtained using an RK4 integrator
with a timestep of $\Delta t = 0.0001$.  Based on the stability results of
Section~\ref{stab_results}, we limit our convergence study to the third-order
schemes, which have the highest potential to improve performance.  We see in
Table~\ref{t:ab34conv} that an estimated order of convergence of about 3 is
obtained in all cases, as expected.

\begin{table}[h] \begin{center}
  \begin{tabular}{ccccc}
       \toprule
       $\text{SR}$ &   err, $\Delta t = 0.005$ & err, $\Delta t = 0.001$ &  err, $\Delta t = 0.0005$ & EOC   \\
       \midrule
     2 &  1.93E-06 & 1.79E-08 & 2.44E-09 & 2.900 \\
     3 &  1.73E-06 & 1.48E-08 & 1.87E-09 & 2.965 \\
     4 &  1.68E-06 & 1.40E-08 & 1.77E-09 & 2.977 \\
     5 &  1.67E-06 & 1.38E-08 & 1.73E-09 & 2.984 \\
     6 &  1.66E-06 & 1.37E-08 & 1.71E-09 & 2.986 \\  
     \bottomrule
 \end{tabular}
 \end{center}
  \caption{Convergence results for third-order MRAB integrators.}\label{t:thirdconv}
\end{table}

\begin{table}[h] \begin{center}
  \begin{tabular}{ccccc}
       \toprule
       $\text{SR}$ &   err, $\Delta t = 0.005$ & err, $\Delta t = 0.001$ &  err, $\Delta t = 0.0005$ & EOC   \\
       \midrule
     2 &  3.56E-06 & 3.42E-08 & 4.53E-09 & 2.894 \\ 
     3 &  3.29E-06 & 2.86E-08 & 3.63E-09 & 2.956 \\ 
     4 &  3.23E-06 & 2.73E-08 & 3.44E-09 & 2.971 \\
     5 &  3.22E-06 & 2.68E-08 & 3.37E-09 & 2.979 \\
     6 &  3.21E-06 & 2.66E-08 & 3.34E-09 & 2.982 \\  
     \bottomrule
 \end{tabular}
 \end{center}
  \caption{Convergence results for third-order MRAB integrators with new extended-history scheme.} \label{t:ab34conv}
\end{table}

\section{MRAB Performance}\label{perf}

\subsection{Performance Model}\label{perf_model}

Our next task is to characterize the benefit in performance when using these
multi-rate integrators.  What expectations should we have in terms of the
reduction in computational work required for these new schemes compared to
single-rate schemes?

To answer this question, we will develop a performance model of how we would
expect an MRAB integrator of a certain step ratio to perform in comparison to
the same RK4 integrator used as a baseline for the accuracy and stability
tests.  In doing so, we make a few assumptions:  we assume that right-hand side
evaluations make up the bulk of the computational cost of the Navier-Stokes
simulation, and we also assume that all simulations will be performed at or
near the maximum stable timestep of a given integrator.

We note in forming this model that in a given timestep for an RK4 integrator,
four right-hand side evaluations must be performed, whereas in the AB schemes,
only one evaluation is needed.  Thus, to reach time $t$, we can write
expressions for the total number of right-hand side evaluations each integrator
will need to perform:
\begin{align}
N_{\text{RHS},{\text{RK4}}} = 4\frac{t}{\Delta t_{\text{RK4}}}(n_{f} + n_{s}) \label{eq:perf_RK4} \\
N_{\text{RHS},{\text{AB}}} = \frac{t}{\Delta t_{\text{RK4}}r_{\text{RK4}}}((\text{SR})n_{f} + n_{s}). \label{eq:perf_AB}
\end{align}
In \eqref{eq:perf_RK4} and \eqref{eq:perf_AB}, $n_{f}$ is the number of points
on the fast grid, and $n_{s}$ is the number of points on the slow grid
(assuming a two-grid configuration).  Furthermore, $\text{SR}$ is the step
ratio of the AB integrator (the number of micro-steps taken on the fast grid
per macro-step on the slow grid --- for an SRAB integrator, $\text{SR} = 1$),
$\Delta t_{\text{RK4}}$ is the timestep taken by the RK4 scheme, and $t$ is the
end time to be reached.  Recall that the factor $r_{\text{RK4}}$ is defined in
Section \ref{stab_results} as 
\begin{align}
r_{\text{RK4}} = \Delta t/\Delta t_{\text{RK4}}. \notag
\end{align}
This factor relates to the timestep restrictions of AB integrators relative to
RK4.  These restrictions, as we have seen in the analyses of
Section~\ref{stab_results}, are imposed by stability constraints.  Namely, the
maximum stable timestep of the RK4 integrator applied to a certain case will be
larger than that of an AB integrator, and therefore more timesteps must be
taken with an AB integrator to reach time $t$ than the RK4 integrator takes to
reach time $t$. 

Given the assumptions of our performance model, calculating the percent
reduction in right-hand side evaluations therefore provides a theoretical
estimate of computational work saved when using a given AB integrator.  We
define the resulting speedup as
\begin{align}
\text{SU} = \frac{N_{\text{RHS},{\text{RK4}}}}{N_{\text{RHS},{\text{AB}}}}. \label{eq:speedup}
\end{align}
Using this definition, any values of SU over 1 indicate a profitable
integrator, whereas values of 1 or less indicate equal or decreased
performance.  As an example of this performance model in execution, we can
compose a model for the case discussed in Section~\ref{case} modeling the
cylinder in crossflow, for a number of step ratios, and for both third and
fourth order.  We use the grid data presented in Section~\ref{case} and the
$r_{\text{RK4}}$ values given in Section~\ref{stab_results}, and for simplicity
we take $t = \Delta t_{\text{RK4}}$.  For these results and for all of the
results that follow, we will use the more potentially profitable
extended-history (AB34 and AB45) integrators.

\begin{table}[h] \begin{center}
  \begin{tabular}{cccc}
       \toprule
       Integrator &                    Total RHS Evals       & \% Red. from RK4      &   SU       \\
       \midrule
     RK4    &                                  34244 &             0.00 &                       1.00  \\ 
     SRAB ($\text{SR} = 1$)&                    33023 &             3.57 &                       1.04  \\
     MRAB ($\text{SR} = 2$) &                   25846 &             24.52 &                     1.32 \\
     MRAB ($\text{SR} = 3$) &                   23529 &             31.29 &                     1.46 \\
     MRAB ($\text{SR} = 4$) &                   22310 &             34.85 &                     1.53 \\                
     MRAB ($\text{SR} = 5$) &                   21581 &             36.98 &                     1.59 \\
     MRAB ($\text{SR} = 6$) &                   25274 &             26.19 &                     1.35 \\
     \bottomrule                
 \end{tabular}
 \end{center}
  \caption{Evaluating RHS costs and theoretical performance for third-order MRAB integrators (using extended history) compared to RK4.}  \label{t:perfmodelthird}
\end{table}

\begin{table}[h] \begin{center}
  \begin{tabular}{cccc}
       \toprule
       Integrator &                               Total RHS Evals       & \% Red. from RK4      &   SU       \\
       \midrule
     RK4  &                                             34244 &             0.00 &                       1.00  \\
     SRAB ($\text{SR} = 1$)&                              48154 &             $-$40.62 &                    0.71  \\
     MRAB ($\text{SR} = 2$) &                            37694 &             $-$10.08 &                      0.91 \\ 
     MRAB ($\text{SR} = 3$) &                            34204 &             0.12 &                       1.00 \\
     MRAB ($\text{SR} = 4$) &                            32459 &             5.21 &                     1.05 \\                
     MRAB ($\text{SR} = 5$) &                            33756 &             1.42 &                     1.01 \\ 
     MRAB ($\text{SR} = 6$) &                            39608 &             $-$15.66 &                  0.86 \\ 
     \bottomrule 
 \end{tabular}
 \end{center}
  \caption{Evaluating RHS costs and theoretical performance for fourth-order MRAB integrators (using extended history) compared to RK4.}  \label{t:perfmodelfourth}
\end{table}

Table~\ref{t:perfmodelthird} shows that we can expect speedup for all
third-order MRAB integrators using the extended history scheme for this
specific case --- however, based on Table~\ref{t:perfmodelfourth}, we speculate
that fourth-order MRAB integrators using the extended histories will largely
fail to be profitable for this specific case, though the fourth-order schemes
do reach a break-even point at $\text{SR}=3$, attaining very minimal
performance benefit at $\text{SR}=4$ (around 5\% reduction in right-hand side
evaluations), before decreasing again due to the stability ceiling observed in
Section~\ref{stab_results}.

\subsection{Sequential Performance} \label{perf_seq}

Having developed a performance model in the previous subsection and thus
obtained expectations for how the integrators should perform, we now develop
tests to verify this model.  We gauge performance by using a given integrator
to march the viscous flow over a cylinder case for 10 NDTU (again, starting from a
steady state solution) at the maximum stable timestep of the integrator,
measuring end-to-end wall clock time of the Navier-Stokes solver, time spent
calculating right-hand sides, time spent performing the necessary
interpolations between grids and time spent performing spatial operator-related
tasks, a subset of the right-hand side calculation.  We limit our performance
testing here to third-order MRAB integrators with the extended history scheme
only, as Section \ref{perf_model} shows these to be the most promising for this
case.  The percentages reported for each step ratio document the percent
reduction in time compared to that of the RK4 integrator, and the predicted
percentage is that which results from the performance model (based purely on
right-hand side evaluations as given by \eqref{eq:perf_RK4} and
\eqref{eq:perf_AB}).

\begin{table}[h] \begin{center}
  \begin{tabular}{cccccccccc}
       \toprule
       Int. &                                              R [s] & \%    &   O [s] &  \%  & I [s] &  \%  & E [s] &  \%   & \% Pred.    \\
       \midrule
     RK4  &                                                          29.9 &    &                   18.4 &                &         0.623   &                   &         38.7 &    & \\ 
    ($\text{SR}1$)&                             28.29 &    5.6 &          17.44 &         5.2 &         0.273  &        56.2 &       41.73                & $-$8.0   & 3.6 \\ 
    ($\text{SR}2$) &                            22.76 &    23.9 &         14.07 &        23.5 &          0.214 &      65.7 &       30.14                & 22.0  & 24.5 \\
    ($\text{SR}3$) &                            21.60 &    27.7 &          13.33 &       27.5 &           0.201 &     67.7 &       28.90               & 25.2   & 31.3 \\ 
    ($\text{SR}4$) &                            20.91 &    30.0 &        12.91 &         29.8 &         0.194 &       68.9 &       27.50               & 28.9   & 34.9 \\                 
    ($\text{SR}5$) &                            18.13 &    39.3 &         11.17 &        39.2 &          0.163 &      73.8 &       23.60                & 38.9  & 37.0 \\
    ($\text{SR}6$) &                            20.27 &    32.2 &         12.51 &        32.0 &          0.182 &      70.8 &       24.13               & 37.6   & 26.2 \\
    \bottomrule                          
 \end{tabular}
 \end{center}
  \caption{Sequential times for third-order MRAB integrators (with extended history).  \textbf{R} = Right-hand Side Evaluation Time, \textbf{O} = Operator Time, \textbf{I} = Interpolation Time, and \textbf{E} = End-to-End Time.}  \label{t:perfthird}
\end{table}

We see in Table~\ref{t:perfthird} that the measured times given for our new
SRAB and MRAB integrators are all lower than those of RK4, with one exception
(the SRAB integrator has a slightly larger end-to-end time).  Furthermore, we
note in the case of sequential performance profiling that the observed
interpolation times are low as we would expect (when not running the solver in
parallel, this amounts to copying data between buffers) --- we will see that
the time spent interpolating becomes more significant as we transition to a
discussion of parallel runs.  Regarding the percentages, while we see good
agreement with the overall trend of the model of Section~\ref{perf_model}, we
see that a number of the percentage reductions observed are slightly lower than
those predicted by the model.

\subsection{Parallel Performance} \label{perf_parallel_small}

In this section, we perform tests similar to those of Section~\ref{perf_seq},
with a specific focus on the distributed-memory parallelization of the new
integrators, in order to answer the question of how communication and
interpolation between grids affects performance.  The procedure with these
tests is the same as in Section \ref{perf_seq}, but an important distinction to
note is that we document end-to-end wall clock time, but the times
reported for right-hand side calculations, operator tasks, and interpolation
are accumulated inclusive times (the total time spent by \emph{all} processors
performing these specific tasks).  In the tables that follow, we report times
for third-order MRAB integrators with extended history for various step ratios.
Bold values indicate lower time spent with a certain task than that of the RK4
integrator.

\begin{table}[h!]  \begin{adjustwidth}{-.5in}{-.5in}  
 \begin{center}
  \begin{tabular}{cccccccccc}
       \toprule
       Proc &  RK4  & $\text{SR}=1$ &                    \% &                   $\text{SR}=2$ &                       \% &                   $\text{SR}=4$ &                   \% &                    $\text{SR}=6$ &                \%        \\
       \midrule
       2 &       23.52 &  25.21 & $-$7.2 & 33.82                 &      $-$43.8       &  25.70          & $-$9.3              & 24.28             & $-$3.2  \\ 
       4 &       17.12 &  19.49 & $-$13.8 & \textbf{10.04}  & \textbf{41.4} & \textbf{8.80} & \textbf{48.6} & \textbf{10.05} & \textbf{41.3} \\
       8 &       8.37 &  9.36 & $-$11.9 & \textbf{6.41}       & \textbf{23.4} & \textbf{6.64} & \textbf{20.6} & \textbf{7.04} & \textbf{15.8} \\
       16 &     4.98 &  5.50 & $-$10.4 & \textbf{3.93}       & \textbf{21.1} & \textbf{3.91} & \textbf{21.5} & \textbf{4.50} & \textbf{9.6} \\        
       \bottomrule
 \end{tabular}
 \end{center}
  \caption{End-to-end timings and percentage reductions for various processor counts.  Bold values indicate lower times than those of the RK4 integrator.}  \label{t:perfendtoend}
  \end{adjustwidth}
\end{table}

\begin{table}[h!] \begin{adjustwidth}{-.5in}{-.5in}  
 \begin{center}
  \begin{tabular}{cccccccccc}
       \toprule
       Proc &  RK4  & $\text{SR}=1$ &                    \% &                   $\text{SR}=2$ &                       \% &                             $\text{SR}=4$ &                   \% &                    $\text{SR}=6$ &                \%        \\
       \midrule
       2 &      31.17 &  \textbf{29.82} & \textbf{4.3} & \textbf{23.67} & \textbf{24.1} & \textbf{20.91} & \textbf{32.9} & \textbf{20.66} & \textbf{33.7} \\ 
       4 &      35.26 &  \textbf{33.40} & \textbf{5.3} & \textbf{27.07}  & \textbf{23.2} & \textbf{24.12} & \textbf{31.6} & \textbf{26.62} & \textbf{24.5} \\
       8 &     44.28  &  \textbf{41.97} & \textbf{5.2} & \textbf{32.97} & \textbf{25.5} & \textbf{28.06} & \textbf{36.6} & \textbf{31.74} & \textbf{28.3} \\ 
       16 &    53.36 &  \textbf{48.82} & \textbf{8.5} & \textbf{39.81} & \textbf{25.4} & \textbf{34.30} & \textbf{35.7} & \textbf{38.97} & \textbf{27.0} \\         
       \bottomrule
 \end{tabular}
 \end{center}
  \caption{RHS timings and percentage reductions for various processor counts.  Bold values indicate lower times than those of the RK4 integrator.}  \label{t:perfrhs}
  \end{adjustwidth}
\end{table}

\begin{table}[h!]  \begin{adjustwidth}{-.5in}{-.5in}  
 \begin{center}
  \begin{tabular}{cccccccccc}
       \toprule
       Proc &  RK4  & $\text{SR}=1$ &                    \% &                   $\text{SR}=2$ &                       \% &                   $\text{SR}=4$ &                   \% &                    $\text{SR}=6$ &                \%        \\
       \midrule
       2 &      19.17 &  \textbf{18.46} & \textbf{3.7} & \textbf{14.61}  & \textbf{23.8} & \textbf{12.86} & \textbf{32.9} & \textbf{12.65} & \textbf{34.0} \\
       4 &      20.97 &  \textbf{20.12} & \textbf{4.1} & \textbf{16.37}  & \textbf{21.9} & \textbf{14.63} & \textbf{30.2} & \textbf{15.78} & \textbf{24.7} \\
       8 &      24.97 &  \textbf{23.83} & \textbf{4.6} & \textbf{18.88}  & \textbf{24.4} & \textbf{16.29} & \textbf{34.8} & \textbf{18.60} & \textbf{25.5} \\
       16 &    28.52 &  \textbf{27.00} & \textbf{5.3} & \textbf{21.71}  & \textbf{23.9} & \textbf{18.29} & \textbf{35.9} & \textbf{21.42} & \textbf{24.9}\\        
       \bottomrule
 \end{tabular}
 \end{center}
  \caption{Operator timings and percentage reductions for various processor counts.  Bold values indicate lower times than those of the RK4 integrator.}  \label{t:perfoperator}
    \end{adjustwidth} 
\end{table}

\begin{table}[h!] \begin{adjustwidth}{-.5in}{-.5in}  
 \begin{center}
  \begin{tabular}{cccccccccc}
       \toprule
       Proc &  RK4  & $\text{SR}=1$ &                    \% &                   $\text{SR}=2$ &                       \% &                   $\text{SR}=4$ &                   \% &                    $\text{SR}=6$ &                \%        \\
       \midrule
       2 &      5.55 &  \textbf{3.86}     & \textbf{30.5} &  29.61          & $-$433.5          & 21.92            & $-$295.0         & 22.04            & $-$297.1\\ 
       4 &      24.54 &  \textbf{21.02} & \textbf{14.3} & \textbf{3.73} & \textbf{84.8} & \textbf{2.28} & \textbf{90.7} & \textbf{5.03} & \textbf{79.5} \\ 
       8 &      13.27 &  \textbf{10.05} & \textbf{24.3} & \textbf{1.46}  & \textbf{89.0} & \textbf{6.60} & \textbf{50.3} & \textbf{11.31} & \textbf{14.8} \\ 
       16 &    15.45 &  \textbf{12.12} & \textbf{21.6} & \textbf{1.93}  & \textbf{87.5} & \textbf{8.53} & \textbf{44.8} & \textbf{14.54} & \textbf{5.9} \\ 
       \bottomrule       
 \end{tabular}
 \end{center}
  \caption{Interpolation timings and percentage reductions for various processor counts.  Bold values indicate lower times than those of the RK4 integrator.}  \label{t:perfinterp}
    \end{adjustwidth}
\end{table}

Generally, what we see here is first and foremost that the time spent on
right-hand side and operator-related tasks (Table~\ref{t:perfrhs} and
Table~\ref{t:perfoperator}, respectively) is always lowered by the use of
multi-rate methods, in certain cases by over 30\%, in line with the model of
Section~\ref{perf_model}).  Furthermore, the end-to-end wall clock times
(Table~\ref{t:perfendtoend}) largely show that the reduction in right-hand side
evaluations indeed leads to end-to-end speedup for a number of step ratios and
processor counts.  Regarding the SRAB scheme ($\text{SR}=1$), we see that the
times spent in the interpolation, operator, and right-hand side tasks are all
lower than those of RK4, yet the end-to-end times are very slightly higher for
all core counts.  The reason for this is the introduction of higher overhead
from elsewhere in the program due to the requirement that more timesteps be run
to reach the same end time.

For tests using 2 processors, we see a large amount of time spent in
interpolation subroutines for all multi-rate integrators. This is due to the
fact that the usage of multi-rate integration in the overset sense
induces load imbalance within the application, given that in a given
macro-timestep, we will be evaluating more right-hand sides on certain grids
than on others.  This causes long processor wait times (accounted for, in our
case, as interpolation time) due to suboptimal work distribution and motivates
a change to the existing decomposition of the problem, which simply distributes
processors to grids based on the ratio of that grid's number of points to the
total number of points in the simulation.  Rescaling this decomposition based
on the multi-rate step ratio being used (a direct indication of how many
right-hand side evaluations per macro-timestep a given processor is responsible
for) allows us to produce improved performance results at higher processor
counts, reducing the wait times spent in grid-to-grid communication.  With only
two cores, however, we are unable to rescale the decomposition to improve
the load imbalance produced by multi-rate. Further effects of inter-grid
communications are shown in the next section, focusing on large-scale results.

\subsection{Large-Scale Performance} \label{perf_parallel_large}

\begin{figure}[h!]
\centering
  \includegraphics[width=0.8\linewidth]{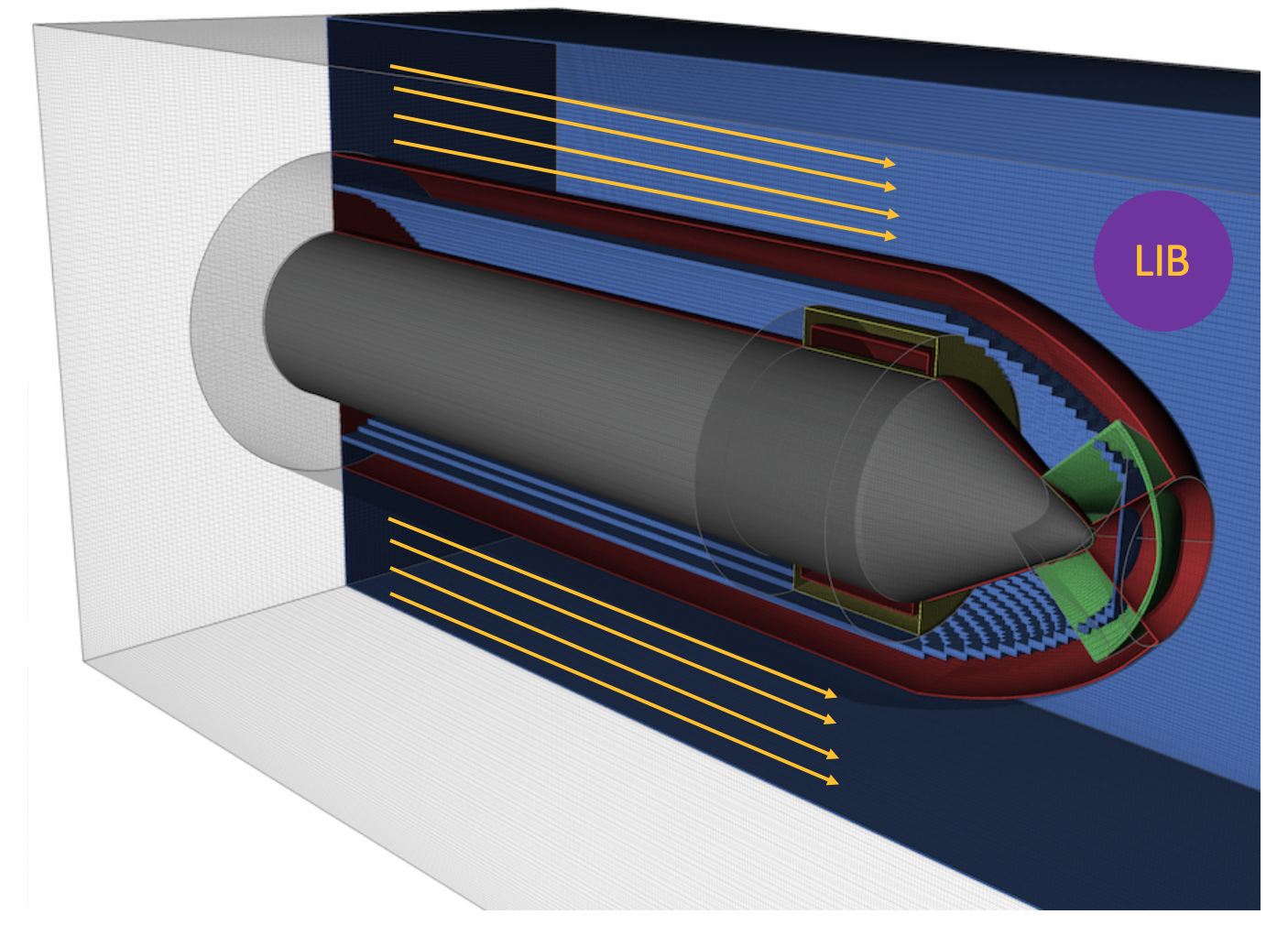}
  \caption{Grids for large-scale case: LIB ignition in underexpanded jet.}
  \label{fig:y4grids}
\end{figure}

\begin{figure}[h!]
 \centering
  \includegraphics[width=0.8\linewidth]{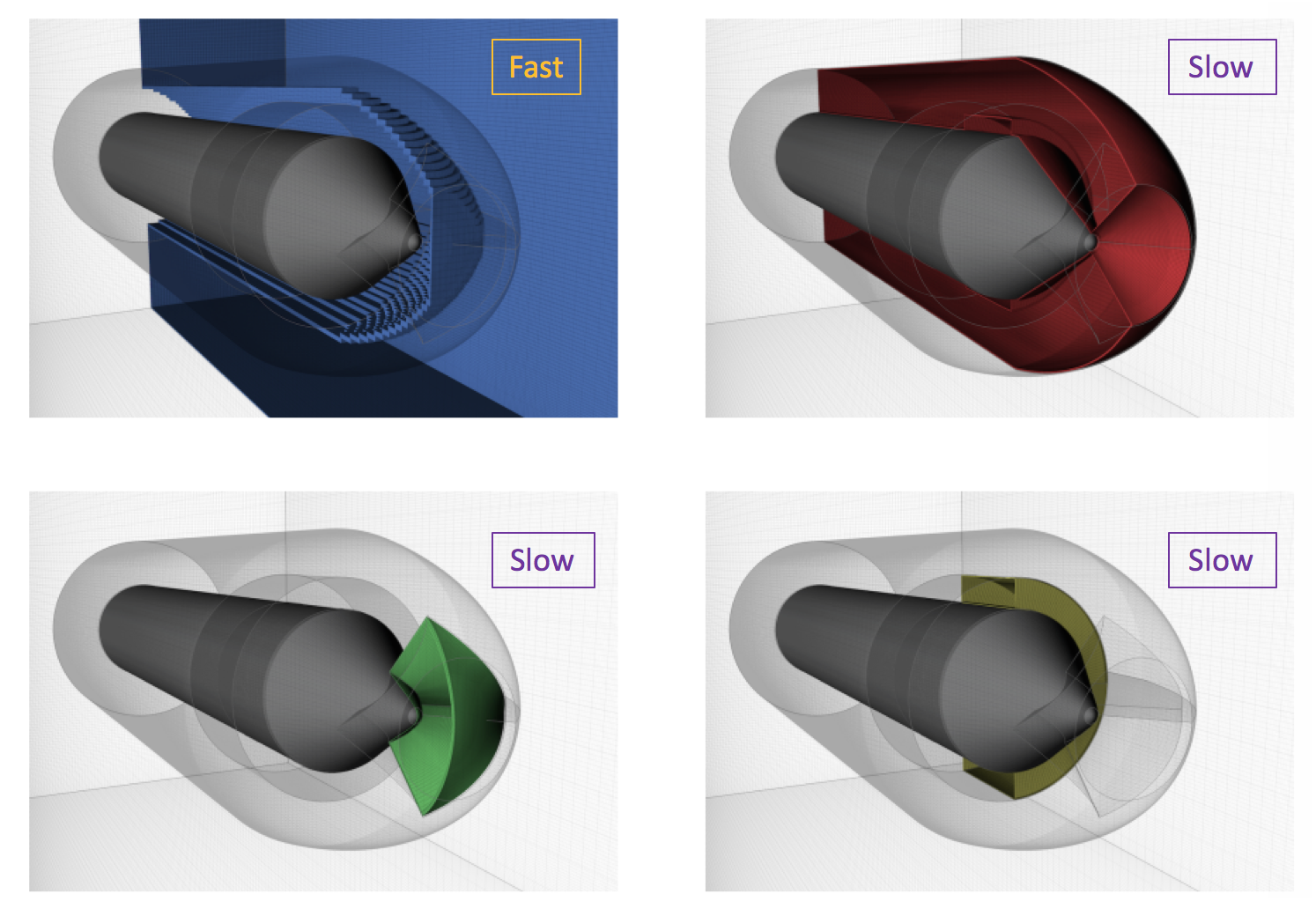}
  \caption{Grids for large-scale case: assigning rates.  Grid 1 (upper left) is assigned the "fast" rate, while Grid 2 (lower left), Grid 3 (upper right), and Grid 4 (lower right) are assigned the "slow" rate.}
  \label{y4rates}
\label{fig:y4ins}
\end{figure}

We now move to a different example to investigate the performance of the
integrators at larger scale.  The grids for this new case, modeling ignition
via laser-induced breakdown (LIB) of an underexpanded jet, are shown in
Figure~\ref{fig:y4grids}, and are numerically described in
Table~\ref{t:perfmodelgrids}.  Note the percentage of total points in Grid 1
--- this will be our fast grid, while the remaining 3 grids will be the slow
grids (see Figure \ref{fig:y4ins}).

\begin{table}[h] \begin{center}
  \begin{tabular}{cccc}
       \toprule
       Grid &                        Grid Type       & No. of Points      &   \% of Total       \\
       \midrule
     1 &                              Cartesian &               12,582,912 &                   47.7  \\
     2 &                              Cylindrical &              339,521 &                        1.3       \\
     3 &                              Cylindrical &               11,419,614 &                   43.3  \\        
     4 &                              Cylindrical &               2,060,250 &                     7.8  \\
     \bottomrule                   
 \end{tabular}
 \end{center}
  \caption{Description of grids for the large-scale case.}  \label{t:perfmodelgrids}
\end{table}


Using the performance model outlined in \eqref{eq:perf_RK4},
\eqref{eq:perf_AB}, and \eqref{eq:speedup}, and assuming $r_{\text{RK4}}$
values for this case similar to the viscous flow over a cylinder case, we can
calculate expected performance.  We test step ratios from 1 to 5, as we have
seen in Section \ref{stab_results} that this is the range of step ratios for
which we maintain perfect efficiency ($r_{\text{SRAB}} \approx \text{SR}$).

\begin{table}[h] \begin{center}
  \begin{tabular}{cccc}
       \toprule
       Integrator &                                           Total RHS Evals       & \% Red. from RK4      &   SU       \\
       \midrule
     RK4 &                                                       105,609,188 &               0.00 &                          1.00  \\ 
     SRAB &                                                    101,939,371 &               3.50 &                          1.04       \\ 
     MRAB ($\text{SR}=2$) &                                      75,261,021 &               26.5 &                          1.36 \\
     MRAB ($\text{SR}=3$) &                                      66,539,511 &               37.0 &                          1.59 \\ 
     MRAB ($\text{SR}=4$) &                                      62,011,632 &               41.3 &                          1.70 \\                     
     MRAB ($\text{SR}=5$) &                                      59,299,803 &               43.8 &                          1.78 \\      
     \bottomrule            
 \end{tabular}
 \end{center}
  \caption{Evaluating RHS costs for third-order MRAB integrators (using extended-history scheme) compared to RK4 for the new large-scale case.}  \label{t:perfmodellarge}
\end{table}

The results in Table~\ref{t:perfmodellarge} suggest that we should expect
moderate speedup when using a third-order multi-rate scheme with $\text{SR}=2$,
with the speedup increasing as we increase the step ratio to 5.  Performing the
integration using these schemes allows us to measure the same times of interest
as in the previous section and compare them with this model.  The results we
present are obtained from a 6,144-core simulation on Stampede2, an
NSF-supported supercomputer at the University of Texas and the Texas Advanced
Computing Center (TACC), and report scaled time spent in various portions of
the Navier-Stokes solver relative to the end-to-end wall clock time of a
Runge-Kutta-driven simulation.

\begin{figure}[h!]
\centering
\includegraphics[width=0.95\textwidth]{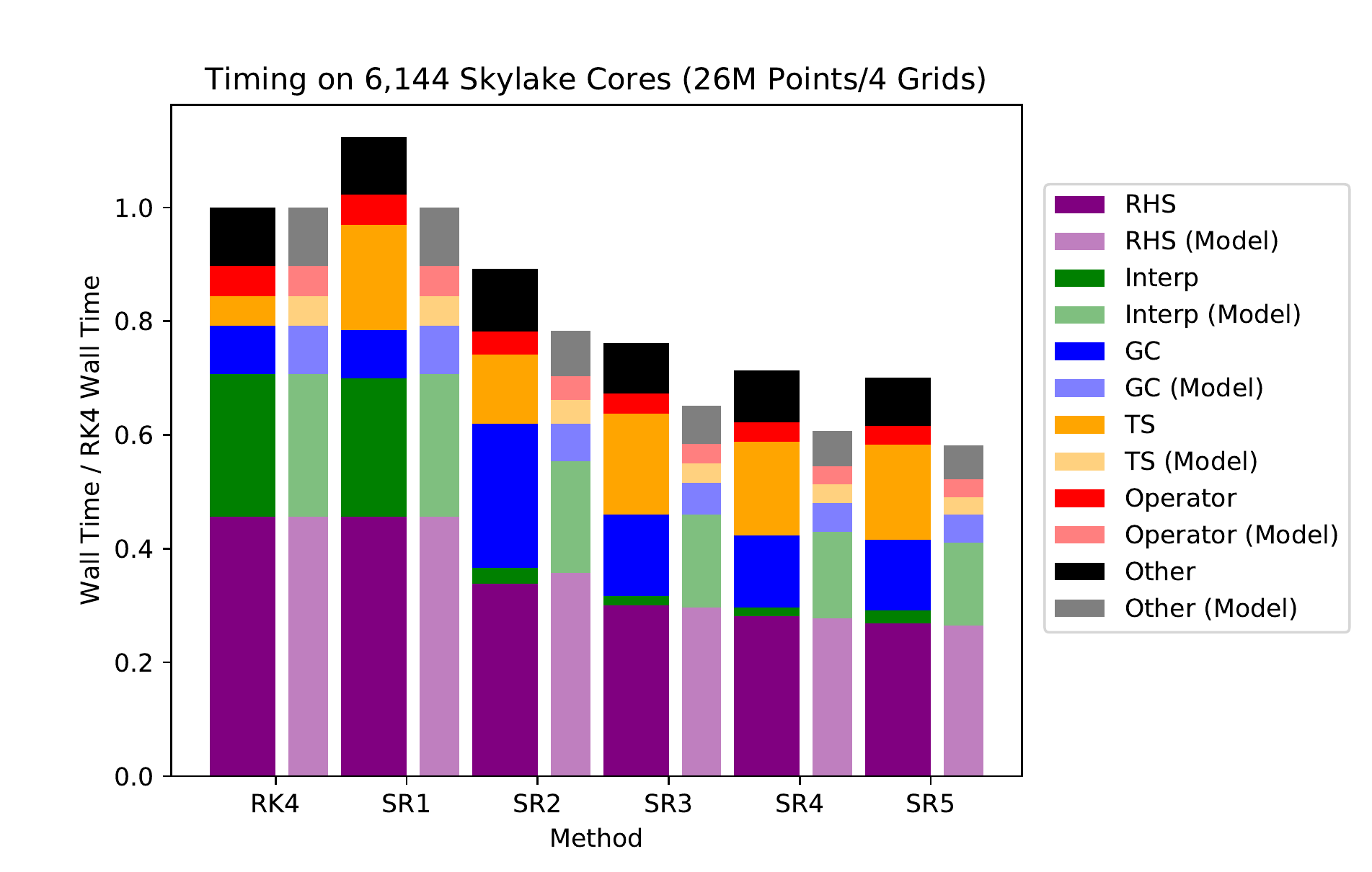}
\caption{Scaled time (wall time / RK4 total wall time) spent in various
components of the flow solver for the large-scale case, comparing RK4
integration and AB34 integration with step ratios of 1-5. Thicker columns
indicate actual performance at a given step ratio, whereas the thinner columns
report model predictions.  "GC" denotes the time spent in ghost-cell exchanges,
and "TS" denotes the time spent in timestep calculation.  We see actual
end-to-end timings above the model prediction for all step ratios, though
actual right-hand-side and operator times are seen to match the model.}
\label{fig:timinggood}
\end{figure}

These results show measured times that correlate well with what we would expect
for reduction in operator and right-hand side costs, and furthermore, we see
large reductions in the time spent in interpolation-related subroutines for all
step ratios --- this is a direct result of implementation of
multi-rate-specific overset interpolation schemes featuring selective
communications and separate send-receive schemes (interleaved with right-hand
side calculation).  In the end, the main time marching loop sees moderate
speedup for all step ratios above 1, but is below the prediction of the
performance model.  Our results show roughly constant time spent in all other
components of the Navier-Stokes solver not explicitly timed, and higher
multi-rate times in two specific portions of the solver, both of which can be
readily explained:

\begin{itemize}
\item \textbf{Timestep calculation.}  While we do scale the frequency of these
	calculations such that they occur the same number of times in RK4 as in
	our integrators, the aforementioned new interpolation scheme --- which
	reduces the amount of wait times in inter-grid communication --- shifts
	load imbalance to these calculations, which feature communication
	between all processors to determine the minimum timestep. 
\item \textbf{Ghost cell exchanges.}  These operations involve communication
	between processors working on the same grid (to exchange ghost cell
	data used to fill boundary stencils).  Rescaling the decomposition the
	solver uses to allocate more processing power to the fast grid and
	account for the inherent load imbalance (as mentioned at the end of
	Section~\ref{perf_parallel_small}) results in more processor boundaries
	that require this operation to be performed.  The question of balancing
	inter-grid communication with intra-grid communication when decomposing
	the problem to avoid load imbalance remains.  We leave this for future
	work.
\end{itemize}

In short, the results for this large-scale case show that while we see
a performance benefit due to decrease in right-hand side evaluations and
implementation of more efficient interpolation routines, the improvement is
below the estimates of the performance model due to increased communication
times induced by load imbalance.

\section{Conclusions and Future Work} \label{conc}
\label{sec:conclusions}

In this paper, we have developed multi-rate Adams-Bashforth integrators to take
advantage of an overset discretization of the Navier-Stokes equations using SBP
spatial operators, SAT boundary conditions, and an SAT penalty-based
interpolation method.  We improve the temporal stability of these integrators
by introducing extended-history Adams-Bashforth schemes with higher maximum
stable timesteps for ODE systems with eigenvalues along the negative real axis.
These new schemes are shown, using a number of numerical tools and examples, to
be stable and accurate, and are also shown to improve performance and reduce
the computational work required to march a given system in time when compared
to both single-rate Adams-Bashforth integrators and a fourth-order Runge-Kutta
integrator.  \revanswera{A simple performance model is also developed (see
Equations \eqref{eq:perf_RK4}, \eqref{eq:perf_AB}, and \eqref{eq:speedup}) that
shows the dependence of expected speedup via the use of multi-rate
Adams-Bashforth integration on both the allowable step ratio and the ratio of
slow points to fast points. In particular, this model shows that performance
benefit is maximized in cases where the ratio of slow points to fast points is
high, and where the fast points have an evolutionary timescale significantly
smaller than that of the slow points.}  In our demonstration of performance
improvement using these new multi-rate schemes, we also identify load balance as a critical
consideration for further reduction in computational time --- this is left for
future work.

\section*{Acknowledgments}
This material is based in part upon work supported by the Department of Energy,
National Nuclear Security Administration, under Award Number DE-NA0002374
and by the National Science Foundation under grant number CCF-1524433.
Computational support has been provided, in part, by the National Science
Foundation XSEDE resources under grant TG-CTS090004.

\bibliography{paper}
\bibliographystyle{abbrvnat}

\appendix

\section{Software and Reproducibility}

Here, we give brief synopses of the software tools used to implement the methods described in this paper.

\begin{itemize}
\item \textbf{PlasComCM} is a Fortran 90 code written to solve the compressible
	Navier-Stokes equations on overset meshes.  PlasComCM is currently
	being used in the University of Illinois' NNSA and DOE-funded PSAAPII
	center, the Center for Exascale Simulation of Plasma-Coupled Combustion
	(XPACC).  For more on XPACC and its work, see
	\url{https://xpacc.illinois.edu}.  For more on PlasComCM, see
	\url{https://bitbucket.org/xpacc/plascomcm}.
\item \textbf{Leap} is a Python package used to describe integration methods
	(including multi-rate integrators) with flexible algorithms via a
	virtual machine, and is capable of describing both implicit and
	explicit time steppers in the form of instructions that can then be
	passed to Dagrt (see below) to generate Fortran or Python code.  Our
	results have been generated using Git revision \texttt{9382dd35} at
	\url{https://github.com/inducer/leap} \cite{Leap2018}.
\item \textbf{Dagrt}, a second Python package, is a DAG-based runtime system
	which can generate Fortran or Python code implementing the integrators
	described by Leap for a given right-hand-side.  In using this tool, and
	Leap, with a host application, the user needs to describe the data
	types to be operated on, along with the right-hand-sides that the host
	application uses, in a short Python driver. Our results have been
	generated using Git revision \texttt{3ccb3479} at
	\url{https://github.com/inducer/dagrt} \cite{Dagrt2018}.
\end{itemize}

\section{Plotting Approximate Stability Regions} \label{approx_stab}
In order to better pose the question of what the maximum stable timestep is,
we must first characterize and plot the approximate stability regions of the
integrator.  These stability regions are plotted in the complex plane, and
enclose a $\lambda \Delta t$ space (where $\lambda$ are the eigenvalues of the
right-hand side being integrated, and $\Delta t$ is the timestep) within which
an integrator is said to be approximately stable.

We approximate the stability region in the complex plane for a given integrator by marching the simple scalar ODE
\begin{align}
\frac{dy}{dt} = ky \notag
\end{align}
in time, with a timestep of $\Delta t = 1$, and with initial condition $y = 1$ at $t = 0$.  By setting
\begin{align}
k = re^{i\theta} + 0.3 \notag
\end{align}
with $0 \leq \theta \leq 2\pi$, we can approximate an integrator's stability
region by varying $r$ for a given value of $\theta$ and marching for 100
timesteps.  For a given $(r,\theta)$ combination, if before 100 steps are taken
the value of $y$ exceeds 2, we deem the method unstable, and $r$ is decreased.
Otherwise, $r$ is increased for that $\theta$ value until the method becomes
unstable.  For all approximate stability regions given, we define the stability
boundary using a tolerance on $r$ of $10^{-12}$, and we use 500 equally spaced
values of $\theta$ in the range $0 \leq \theta \leq 2\pi$.  This procedure is
essentially an application of the boundary locus method described in
\cite{sanz1980some} for ordinary differential equations, which we find to be a
reasonable heuristic for discussion of stability in the context of
single-rate time integration of PDEs using an SBP-SAT
discretization.

\section{\revanswerbb{Demonstrating Spatial Order of Accuracy}} \label{convergence}

\subsection{Spatial Accuracy of Right-Hand Sides}

To spatially discretize the second-derivative terms in the compressible
Navier-Stokes equations, we use two applications of the first-derivative
operator.  A Taylor series expansion analysis of the SBP 4-2 spatial operators
employed in this paper shows that if the corresponding first-derivative
operator is accurate of order $p$ on the boundary and accurate of order $2p$ in
the interior, then the resulting second-derivative operator is accurate of
order $p-1$ on the boundary while retaining order $2p$ in the interior.  To
calculate the global error in two dimensions for a Cartesian domain, we use the expression
\begin{align}
\|\boldsymbol{\epsilon}\|_{\tilde{2}} = \sqrt{\sum\limits_{k=1}^M
\sum\limits_{i=1}^{N_{i,k}}\sum\limits_{j=1}^{N_{j,k}}
(\epsilon_{i,j}^k)^{2}\Delta x_{k} \Delta y_{k}}, \label{eq:l2error}
\end{align}
where $M$ is the number of subdomains, $N_{i,k}$ is the number of points in the
$i$-direction on the $k$-th subdomain, and $N_{j,k}$ is the number of points in
the $j$-direction on the $k$-th subdomain. An analysis of this expression given
the expectations for boundary and interior orders of accuracy shows that the
global order we should expect as measured in this norm is $p+1/2$ for the
first-derivative operator, and $p-1/2$ for the second derivative operator.
This result is stated in \citep{svard2014review} while also noting that
solution stability guarantees higher global order in practice when the
solutions are used to calculate error
\citep{gustafsson1975convergence,gustafsson1981convergence}, and when the Kreiss
condition \citep{gustafsson1995time} is satisfied.

To demonstrate this design order, we will first calculate the error in
right-hand sides for prescribed solutions for which the analytical
right-hand sides of the Navier-Stokes equations are known.  The compressible
Navier-Stokes equations (without source terms or heat fluxes) are given by 
\begin{align}
\frac{\partial \rho}{\partial t} + \frac{\partial}{\partial x_{j}} \rho u_{j} &= 0 \notag \\
\frac{\partial \rho u_{i}}{\partial t} + \frac{\partial}{\partial x_{j}} (\rho u_{i} u_{j} + p \delta_{ij} - \tau_{ij}) &= 0 \notag \\
\frac{\partial \rho E}{\partial t} + \frac{\partial}{\partial x_{j}} ((\rho E + p)u_{j} - u_{i}\tau_{ij}) &= 0. \notag
\end{align}
Here, we have viscous stress terms given by
\begin{align}
\tau_{ij} = \frac{\mu}{\text{Re}}\left(\frac{\partial u_{i}}{\partial x_{j}} + \frac{\partial u_{j}}{\partial x_{i}}\right) + \frac{\lambda}{\text{Re}}\frac{\partial u_{k}}{\partial x_{k}}\delta_{ij}. \notag
\end{align}
Setting an initial condition of 
\begin{align}
\rho &= 1 \notag \\
\rho u &= \sin(2\pi (x - 0.5)) \notag \\
\rho v &= 0 \notag \\
\rho E &= \frac{p}{1 - \gamma} + \frac{1}{2} \rho (u^2 + v^2). \notag
\end{align}
we expect the initial time derivative of the density to be equal to
\begin{align}
\frac{\partial \bar{\rho}}{\partial t} = 2\pi \cos(2\pi (x - 0.5)). \notag
\end{align}
We calculate the error in the time derivative of the density via Equation~\ref{eq:l2error}, where
 \begin{align}
 \epsilon_{i,j}^k = \left(\frac{\partial \rho}{\partial t}\right)_{i,j}^k - \left(\frac{\partial \bar{\rho}}{\partial t}\right)_{i,j}^k. \notag
\end{align}
The domain upon which we test these manufactured right-hand sides uses a
box-in-box grid configuration.  All boundaries on the outer grid are periodic.
The outer grid has dimensions $[-1,1] \times [-1,1]$ and the inner grid has
dimensions $[-0.5,0.5] \times [-0.5,0.5]$.  The outer grid is discretized using
$N \times N$ uniformly distributed points, whereas the inner grid is
discretized using $N/2 \times N/2$ uniformly distributed points.  The
interpolation between grids for this case and all others in this appendix uses
the SAT penalty-based scheme described in Section~\ref{interp}.
Table~\ref{t:mms_rhs1} shows that we obtain a global order of accuracy of
$p+1/2$ for the inviscid problem.
\begin{table}[h] \begin{center}
  \begin{tabular}{ccc}
       \toprule
       N & log10($\|\boldsymbol{\epsilon}\|_{2}$) & EOC   \\
       \midrule
     101 & -2.0518 & - \\ 
     201 & -2.7834 & 2.448 \\ 
     401 & -3.5275 & 2.481 \\
     \bottomrule
 \end{tabular}
 \end{center}
  \caption{Convergence results using the method of manufactured solutions - inviscid-only RHS with SBP4-2 first-derivative operators.}  \label{t:mms_rhs1}
\end{table}

For the viscous problem, we prescribe the following solution, which admits no inviscid (first-derivative) right-hand side components.
\begin{align}
\rho &= 1 \notag \\
\rho u &= \sin(2\pi (y - 0.5)) \notag \\
\rho v &= 0 \notag \\
\rho E &= \frac{p}{1 - \gamma} + \frac{1}{2} \rho (u^2 + v^2) \notag.
\end{align}
With this initial condition, we expect the time derivative of the $x$-momentum component to be equal to 
\begin{align}
\frac{\partial \rho u}{\partial t} = \frac{\partial \tau_{12}}{\partial y} = -4\pi^{2} \sin(2\pi (y - 0.5)). \notag
\end{align}
Performing the same experiment as above on the same spatial domain, and
calculating the error in the time derivative of the $x$-momentum component
rather than in the density (the density right-hand side admits no
second-derivative terms), we see in Table~\ref{t:mms_rhs2} that here we obtain
a global order of accuracy of roughly $p-1/2$, which is again in line with our
expectation.

\begin{table}[h] \begin{center}
  \begin{tabular}{ccc}
       \toprule
       N & log10($\|\boldsymbol{\epsilon}\|_{2}$) & EOC   \\
       \midrule
     101 & -2.3518 & - \\ 
     201 & -2.7834 & 1.444 \\ 
     401 & -3.2270 & 1.479 \\
     \bottomrule
 \end{tabular}
 \end{center}
  \caption{Convergence results using the method of manufactured solutions for the viscous-only RHS with SBP 4-2 second-derivative operators constructed using repeated first derivatives.}  \label{t:mms_rhs2}
\end{table}

\subsection{Spatial Accuracy of Solutions}

\begin{figure}[h!]
\centering
\begin{subfigure}[b]{0.50\textwidth}
 \centering
  \includegraphics[width=0.9\linewidth]{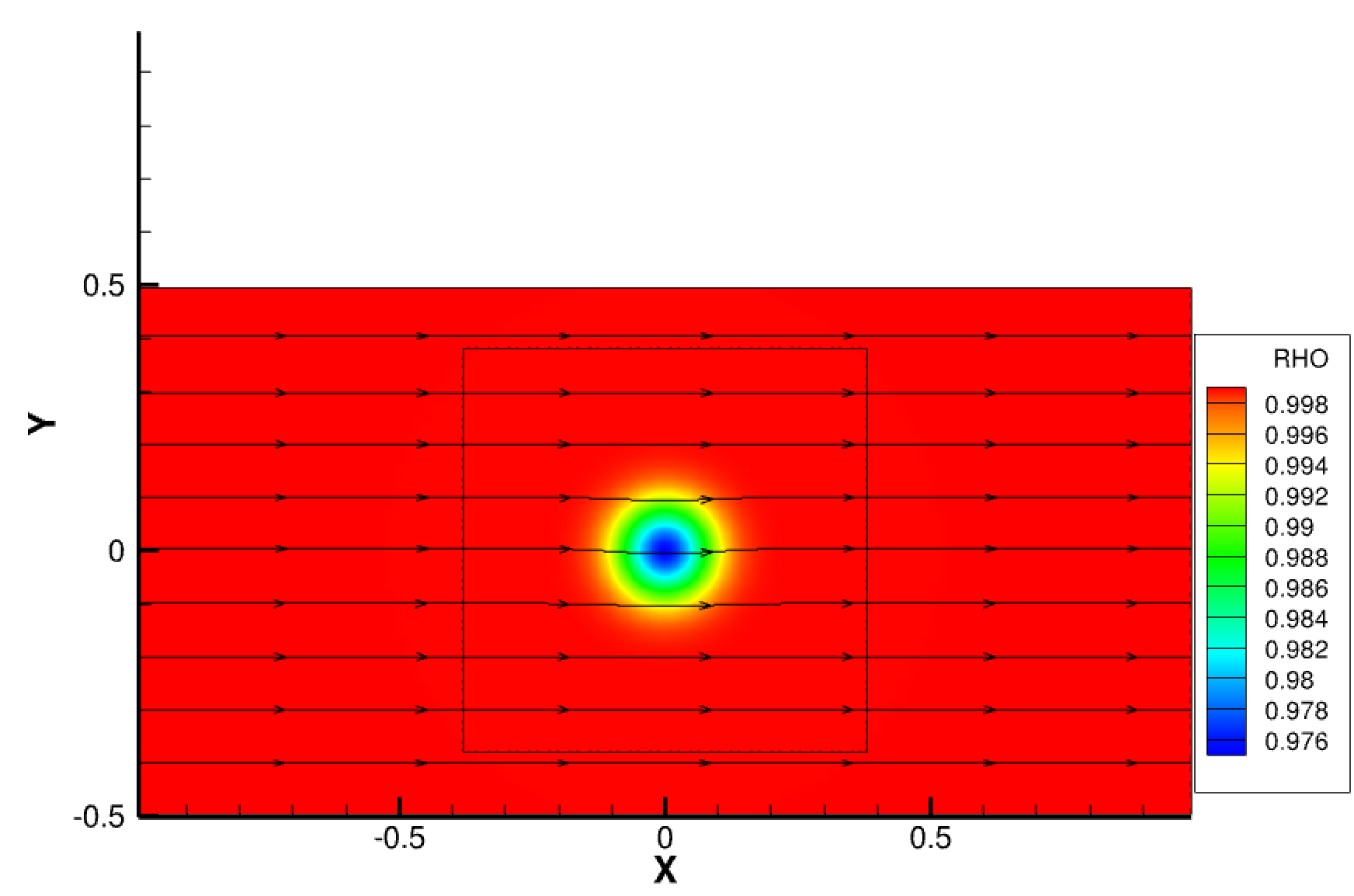}
  \caption{Box-in-box grid configuration for \newline convecting vortex test.}
  \label{fig:inviscid_grids}
\end{subfigure}\begin{subfigure}[b]{0.50\textwidth}
 \centering
  \includegraphics[width=0.9\linewidth]{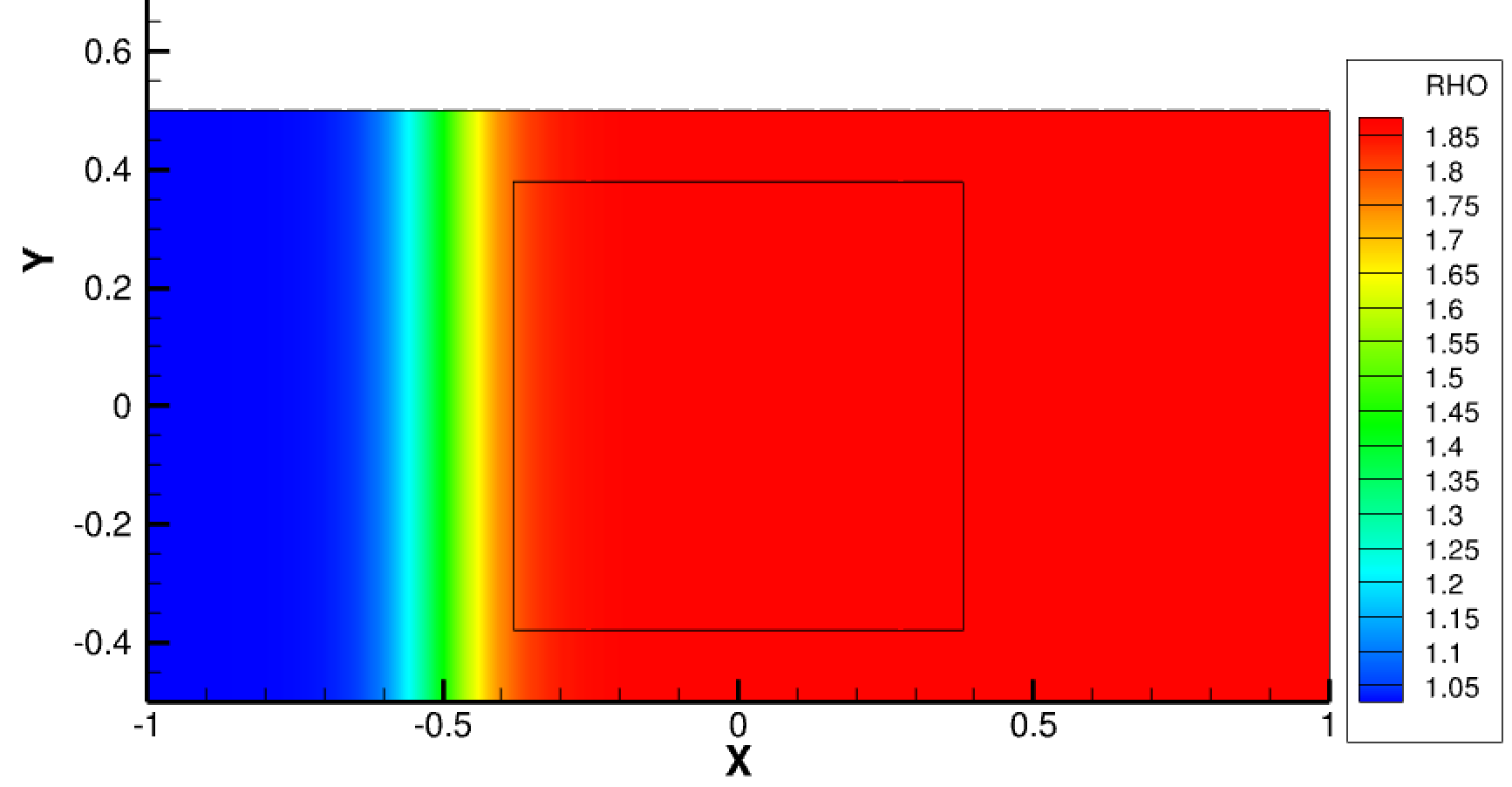}
  \caption{Box-in-box grid configuration for moving shock test.}
  \label{fig:viscous_grids}
\end{subfigure}
\caption{Computational domains for testing of spatial order of accuracy using Navier-Stokes solutions.}
\label{fig:soln_grids}
\end{figure}

We now demonstrate spatial order of accuracy of the SBP 4-2 spatial operators, first using
the same convecting vortex solution of the 2D Euler equations as in~\citep{sharan2018time} on a
box-in-box grid configuration.  The 2D Euler equations are given by
\begin{align}
\frac{\partial \rho}{\partial t} + \frac{\partial}{\partial x_{j}} \rho u_{j} = 0, \notag \\
\frac{\partial \rho u_{i}}{\partial t} + \frac{\partial}{\partial x_{j}} (\rho u_{i} u_{j} + p \delta_{ij}) = 0, \notag \\
\frac{\partial \rho E}{\partial t} + \frac{\partial}{\partial x_{j}} ((\rho E + p)u_{j}) = 0, \notag
\end{align}
and the convecting vortex solution is given by
\begin{align}
\rho &= \left(1 - \frac{\omega^{2}(\gamma - 1)}{8\pi^{2}c_{0}^{2}}e^{1 - \phi^{2}r^{2}}\right)^{\frac{1}{\gamma - 1}}, \notag \hspace{7mm}
v = v_{0} + \frac{\omega}{2\pi} \phi (x - x_{0} - u_{0}t)e^{\frac{1 - \phi^{2}r^{2}}{2}}, \notag \\
u &= u_{0} - \frac{\omega}{2\pi} \phi (y - y_{0} - v_{0}t)e^{\frac{1 - \phi^{2}r^{2}}{2}}, \notag \hspace{5mm}
p = \rho^{\gamma}, \notag \\
\rho E &=  \frac{p}{\gamma - 1} + \frac{1}{2} \rho (u^{2} + v^{2}), \notag \hspace{18mm}
r^{2} = (x - x_{0} - u_{0}t)^{2} + (y - y_{0} - v_{0}t)^{2}, \notag
\end{align}
where $(x_0,y_0)$ is the initial position of the vortex, $(u_{0}, v_{0})$ is
the vortex convective velocity, $\phi$ is a scaling factor, and $\omega$
denotes the non-dimensional circulation.  We set the bulk velocity $(u_{0},
v_{0})$ to be (2,0) and specify periodic boundaries in both the $x$- and $y$-
directions, and we measure the error at $t = 1$.  The outer grid has dimensions
$[-1,1] \times [-0.5,0.5]$ and the inner grid has dimensions $[-0.38, 0.38]
\times [-0.38, 0.38]$, and both grids are discretized using $N \times N$
uniformly distributed points.  For all convergence tests that follow, we
calculate the error in the density to confirm global accuracy.  The
error is calculated via Equation~\ref{eq:l2error}, where
\begin{align}
\epsilon_{i,j}^k = (\rho_{i,j}^k - \bar{\rho}_{i,j}^k), \notag 
\end{align}
and $\bar{\rho}$ is the analytical solution at $t = 1$.  The domain and initial
condition are shown in Figure~\ref{fig:inviscid_grids}.  For SBP 4-2
first-derivative operators, we expect the global accuracy to be of order 3, as
demonstrated in \citep{gustafsson1975convergence}.  Table~\ref{t:mms1} shows
that we attain this order of accuracy for this case.

\begin{table}[h] \begin{center}
  \begin{tabular}{ccc}
       \toprule
       N & log10($\|\boldsymbol{\epsilon}\|_{2}$) & EOC   \\
       \midrule
     50 & -3.0921 & - \\
     100 & -4.0742 & 3.263 \\
     150 & -4.7235 & 3.688 \\
     200 & -5.1681 & 3.558 \\
     \bottomrule  
 \end{tabular}
 \end{center}
  \caption{Convergence results for the 2D Euler equations for the convecting vortex case.}  \label{t:mms1}
\end{table}

We also test a viscous problem to demonstrate design order when
second-derivative terms are present.  The experiment uses the moving shock
problem described by \citep{svard2006order} and \citep{whitham2011linear}, and
is performed on the box-in-box grid configuration shown in
Figure~\ref{fig:viscous_grids}.  As in the inviscid case, the outer grid has
dimensions $[-1,1] \times [-0.5,0.5]$ and the inner grid has dimensions
$[-0.38, 0.38] \times [-0.38, 0.38]$.  In this case, the Reynolds number is
20.8429, the Prandtl number is 0.75, and the shock is moving at $M = 0.5$ relative to the grid,
starting from $x_0 = -0.5$ and propagating to the interior grid.  Once again,
we calculate the error in the density using Equation~\ref{eq:l2error}.  The
errors are measured at time $t = 0.52$.  Table~\ref{t:mms2} shows that we
attain the expected order of accuracy.

\begin{table}[h] \begin{center}
  \begin{tabular}{ccc}
       \toprule
       N & log10($\|\boldsymbol{\epsilon}\|_{2}$) & EOC   \\
       \midrule
     41 &  -2.721 & - \\ 
     81 & -3.613 & 3.014 \\ 
     161 & -4.505 & 2.992 \\
     \bottomrule
 \end{tabular}
 \end{center}
  \caption{Convergence results for the compressible Navier-Stokes equations for the moving shock case.}  \label{t:mms2}
\end{table}

\end{document}